\documentclass[11pt]{amsart}
\usepackage{amsxtra}
\usepackage{amssymb}
\usepackage[mathscr]{eucal}
\theoremstyle{plain}

\newcommand{\cleqn}{\setcounter{equation}{0}}
\newcommand{\clth}{\setcounter{theorem}{0}}
\newcommand {\sectionnew}[1]{\section{#1}\cleqn\clth}

\newtheorem{theorem}{Theorem}[section]
\newtheorem{lemma}[theorem]{Lemma}
\newtheorem{definition-theorem}[theorem]{Definition-Theorem}
\newtheorem{proposition}[theorem]{Proposition}
\newtheorem{corollary}[theorem]{Corollary}
\newtheorem{definition}[theorem]{Definition}
\newtheorem{example}[theorem]{Example}
\newtheorem{remark}[theorem]{Remark}
\newtheorem{conjecture}[theorem]{Conjecture}

\newtheorem*{maintheorem*}{Main Theorem}
\newcommand \bth[1] { \begin{theorem}\label{t#1} }
\newcommand \ble[1] { \begin{lemma}\label{l#1} }

\newcommand \bpr[1] { \begin{proposition}\label{p#1} }
\newcommand \bco[1] { \begin{corollary}\label{c#1} }
\newcommand \bde[1] { \begin{definition}\label{d#1}\rm }
\newcommand \bex[1] { \begin{example}\label{e#1}\rm }
\newcommand \bre[1] { \begin{remark}\label{r#1}\rm }
\newcommand \bcj[1] { \begin{conjecture}\label{j#1}\rm }

\renewcommand {\eth} { \end{theorem} }
\newcommand {\ele} { \end{lemma} }

\newcommand {\epr} { \end{proposition} }
\newcommand {\eco} { \end{corollary} }
\newcommand {\ede} { \end{definition} }
\newcommand {\eex} { \end{example} }
\newcommand {\ere} { \end{remark} }
\newcommand {\ecj} { \end{conjecture} }

\newcommand \thref[1]{Theorem \ref{t#1}}
\newcommand \leref[1]{Lemma \ref{l#1}}
\newcommand \prref[1]{Proposition \ref{p#1}}
\newcommand \coref[1]{Corollary \ref{c#1}}

\newcommand \deref[1]{Definition \ref{d#1}}
\newcommand \exref[1]{Example \ref{e#1}}
\newcommand \reref[1]{Remark \ref{r#1}}

\usepackage{tikz}
\usetikzlibrary{matrix}
\usepackage{graphicx,ctable,booktabs}
\usetikzlibrary{arrows.meta}
\DeclareMathOperator{\up}{up} 
\DeclareMathOperator{\low}{low}
\DeclareMathOperator{\Span}{Span}

\DeclareMathOperator{\Spec}{Spec} \DeclareMathOperator{\Ann}{Ann}

\DeclareMathOperator{\End}{End} 
 \DeclareMathOperator{\Hom}{Hom}

\DeclareMathOperator{\Prim}{Prim}

\newcommand{\mc}{\mathcal}

\newcommand{\Cat}{\operatorname{Cat}}
\newcommand{\id}{\operatorname{id}}

\newcommand{\kk}{\Bbbk}
\def \Cset {{\mathbb C}}

\def \Zset {{\mathbb Z}}

\def \Qset {{\mathbb Q}}

\def \mt {\mapsto}
\def \g {\mathfrak{g}}      
\def \t {\mathfrak{t}} 
\def \al {\alpha} 
\def \be {\beta} 
\def \vpi {\varpi}
\def \la {\lambda}
\def \AA {\mathcal{A}}
\def \BB {\bf{B}}

\begin{document}

\title[Prime Spectra of 2-Categories and 
Categorifications of Richardson Varieties]
{Prime Spectra of Abelian 2-Categories and 
Categorifications of Richardson Varieties}
\dedicatory{To Anthony Joseph on his 75th birthday, with admiration}
\author[Kent Vashaw]{Kent Vashaw}
\address{
Department of Mathematics \\
Louisiana State University \\
Baton Rouge, LA 70803 \\
U.S.A.
}
\email{kvasha1@lsu.edu}
\author[Milen Yakimov]{Milen Yakimov}
\email{yakimov@math.lsu.edu}
\begin{abstract}
We describe a general framework for prime, completely prime, semiprime, and primitive ideals of an abelian 2-category. This provides a noncommutative version of Balmer's 
prime spectrum of a tensor triangulated category. These notions are based on containment conditions in terms of thick subcategories 
of an abelian category and thick ideals of an abelian 2-category. We prove categorical analogs of the main properties of noncommutative prime spectra.
Similar notions, starting with Serre subcategories of an abelian category and Serre ideals of an abelian 2-category, are developed. They are linked to 
Serre prime spectra of $\mathbb{Z}_+$-rings.
As an application, we construct a categorification of the quantized coordinate rings of open Richardson varieties for symmetric Kac--Moody groups, 
by constructing Serre completely prime ideals of monoidal categories of modules of the KLR algebras, and by taking Serre quotients with respect to them.  
\end{abstract}
\maketitle
\sectionnew{Introduction}
\subsection{Noncommutative categorical prime spectra}
Balmer's notion of a prime spectrum of a tensor triangulated category \cite{Balmer1,Balmer2} is a major tool in homological algebra, representation theory, 
algebraic topology and other areas. It is defined for triangulated categories with a symmetric monoidal structure.
As noted in \cite{BKN}, Balmer's construction and results more generally apply to braided monoidal triangulated categories.

The notion of a prime spectrum of a braided monoidal triangulated category is a categorical version of the notion of a prime spectrum of a commutative ring. In the classical 
case of noncommutative rings, there are four different notions of primality \cite{Goodearl1}: prime, completely prime, semiprime and primitive spectra. In this paper we develop 
categorical notions of all of them, and prove analogs of many of their main properties. We do this in the abelian setting. However, instead of simply considering 
abelian monoidal categories, we work with the more general setting of abelian 2-categories. It is necessary to consider this more general setting, because many 
of the monoidal categorifications of noncommutative algebras 
that have been constructed so far are in the setting of 2-categories, rather than monoidal categories, see \cite{Kamnitzer,Lauda1,Mazorchuk1}.
Categorifications via 2-categories are even needed for relatively 
small algebras such as the idempotented version of the quantized universal enveloping algebra of $\mathfrak{sl}_2$; we refer the reader to \cite{Lauda1} for a very 
informative review of this particular categorification.
\subsection{Thick and prime ideals of abelian 2-categories} 
\label{1.2}
A {\em{2-category}} is a category enriched over the category of 1-categories. In other words, a 
2-category $\mc{T}$, has the property that for every two objects $A_1, A_2$ of it, the morphisms $\mc{T}(A_1, A_2)$ form a 1-category 
and satisfy natural 
identity conditions. A 2-category with one object is the same thing as a strict monoidal category.
An {\em{abelian 2-category}} is such a category for which the 1-categories $\mc{T}(A_1, A_2)$ are abelian and the composition bifunctors are 
biexact. We work with small 2-categories, i.e., with 2-categories $\mc{T}$ whose objects form a set and for which all 1-categories 
$\mc{T}(A_1, A_2)$ are small. 
We denote by $\mc{T}_1$ the isomorphism classes of 1-morphisms of $\mc{T}$. 
For two subsets $X, Y \subseteq \mc{T}_1$, denote by 
\begin{align*}
\mbox{$X \circ Y$} \; &\mbox{the set of isomorphism classes of 1-morphims of $\mc{T}$ having representatives of the}
\\
&\mbox{form $f g$ for $f$ and $g$ representing classes in $X$ and $Y$ such that $fg$ is defined}.
\end{align*}

The different versions of prime ideals of abelian 2-categories which we develop are based on the notion of a {\em{thick subcategory of an
abelian category}} and its 2-incarnation, the notion of a {\em{thick ideal of an abelian 2-category}}. Recall that a thick (sometimes called wide) 
subcategory of an abelian category is a nonempty full subcategory which is closed under taking kernels, 
cokernels, and extensions.

A {\em{thick ideal}} $\mc{I}$ of an abelian 2-category $\mc{T}$ is a collection of subcategories $\mc{I}(A_1, A_2)$ of $\mc{T}(A_1, A_2)$ for all objects 
$A_1, A_2$ of $\mc{T}$ such that 
\begin{enumerate}
\item $\mc{I}(A_1, A_2)$ are thick subcategories of the abelian categories $\mc{T}(A_1, A_2)$ and 
\item the composition bifunctors of $\mc{T}$ restrict to bifunctors
\[
\mc{T}(A_2, A_3) \times \mc{I}(A_1, A_2) \to \mc{I}(A_1, A_3) \quad \mbox{and}  \quad
\mc{I}(A_2, A_3) \times \mc{T}(A_1, A_2) \to \mc{I}(A_1, A_3)
\] 
for all objects $A_1, A_2, A_3$ of $\mc{T}$.
\end{enumerate}

We call a proper thick ideal $\mc{P}$ of an abelian 2-category $\mc{T}$
\begin{enumerate}
\item[(p)] {\em{prime}} if for all thick ideals $\mc{I}$ and $\mc{J}$ of $\mc{T}$, $\mc{I}_1 \circ \mc{J}_1 \subseteq \mc{P}_1$ 
implies that either $\mc{I} \subseteq \mc{P}$ or  $\mc{J} \subseteq \mc{P}$,
\item[(sp)] {\em{semiprime}} if it is an intersection of prime ideals,
\item[(cp)] {\em{completely prime}} if, for all $f, g \in \mc{T}_1$, $f \circ g \subseteq \mc{P}_1$ implies that either $f \in \mc{P}_1$ or $g \in \mc{P}_1$.
Note that the set $f \circ g$ is either empty or is a singleton.
\end{enumerate}
We obtain categorical versions of the main properties of prime, semiprime, and completely prime ideals of noncommutative rings. 
In Section \ref{3} it is proved that the following are equivalent for a proper thick ideal $\mc{P}$ of an abelian 2-category $\mc{T}$:
\begin{enumerate}
\item[(p1)] $\mc{P}$ is a prime ideal;
\item[(p2)] If $f, g \in \mc{T}_1$ and $f \circ \mc{T}_1 \circ g \subseteq \mc{P}_1$, then either $f \in \mc{P}_1$ or $g \in \mc{P}_1$;  
\item[(p3)] If $\mc{I}$ and $\mc{J}$ are any thick ideals properly containing $\mc{P}$, then $\mc{I}_1 \circ \mc{J}_1\not \subseteq \mc{P}_1$;
\item[(p4)] If $\mc{I}$ and $\mc{J}$ are any left thick ideals of $\mc{T}$ such that $\mc{I}_1 \circ \mc{J}_1 \subseteq \mc{P}_1$, then either 
$\mc{I} \subseteq \mc{P}$ or $\mc{J} \subseteq \mc{P}$. 
\end{enumerate}
We call the set of such thick ideals of $\mc{T}$ the {\em{prime spectrum}} of $\mc{T}$, to be denoted by $\Spec (\mc{T})$, and define a 
{\em{Zariski type topology}} on it.
For a {\em{multiplicative subset}} $\mc{M}$ of $\mc{T}_1$ (see \deref{mult}) 
and a proper thick ideal $\mc{I}$ of $\mc{T}$ such that $\mc{I}_1 \cap \mc{M} = \varnothing$, we prove
that every {\em{maximal element}} of the set 
\[
X(\mc{M}, \mc{I}) := \{ \mc{K} \; \mbox{a thick ideal of} \; \; \mc{T} \mid \mc{K} \supseteq \mc{I} \; \; \mbox{and} \; \; \mc{K}_1 \cap \mc{M} = \varnothing \}
\]
is a prime ideal of $\mc{T}$. This implies that $\Spec (\mc{T})$ is {\em{non-empty}} for every abelian 2-category.

Categorical versions of simple, noetherian and weakly noetherian noncommutative rings are given in Section \ref{4}. There we prove that for every 
weakly noetherian abelian 2-category $\mc{T}$ and a proper thick ideal $\mc{I}$ of $\mc{T}$, 
there exist finitely many minimal Serre prime ideals over $\mc{I}$ and 
there is a finite list of minimal prime ideals over $\mc{I}$ (possibly with repetition) $\mc{P}^{(1)},..., \mc{P}^{(m)}$ such that the product 
\[
\mc{P}^{(1)}_1 \circ ...\circ \mc{P}^{(m)}_1 \subseteq \mc{I}_1.
\]

In Section \ref{5} we prove a categorical version of the Levitzki--Nagata theorem for semiprime ideals, and furthermore show that 
the following are equivalent for a proper thick ideal $\mc{Q}$ of $\mc{T}$:
\begin{enumerate}
\item[(sp1)] $\mc{Q}$ is semiprime;
\item[(sp2)] If $f \in \mc{T}_1$ and $f \circ \mc{T}_1 \circ f \subseteq \mc{Q}_1$, then $f \in \mc{Q}_1$.
\item[(sp3)] If $\mc{I}$ is any thick ideal of $\mc{T}$ such that $\mc{I}_1 \circ \mc{I}_1 \subseteq \mc{Q}_1$, then $\mc{I} \subseteq \mc{Q}$;
\item[(sp4)] If $\mc{I}$ is any thick ideal properly containing $\mc{Q}$, then $\mc{I}_1 \circ \mc{I}_1\not \subseteq \mc{Q}_1$;
\item[(sp5)] If $\mc{I}$ is any left thick ideal of $\mc{T}$ such that $\mc{I}_1 \circ \mc{I}_1 \subseteq \mc{Q}_1$, then $\mc{I} \subseteq \mc{Q}$.
\end{enumerate}
\subsection{Serre prime ideals of 2-categories and ideals of $\Zset_+$-rings} 
\label{1.3}
Serre subcategories of abelian categories are a particular 
type of thick subcategories. A thick ideal $\mc{I}$ of an abelian 2-category $\mc{T}$ will be called a {\em{Serre ideal}} if $\mc{I}(A_1, A_2)$ is a
Serre subcategory of $\mc{T}(A_1, A_2)$ for all objects $A_1, A_2$ of $\mc{T}$.
For those ideals one can consider the {\em{Serre quotient}} $\mc{T}/\mc{I}$ which is an abelian 2-category under 
a mild condition on $\mc{I}$.

We define a {\em{Serre prime}} (resp. {\em{semiprime, completely prime}}) {\em{ideal}} of an abelian 2-category $\mc{T}$ to be a prime (resp. semiprime, completely prime)
ideal which is a Serre ideal. Section \ref{6} treats in detail these ideals, and proves that they are characterized by similar to (p3)-(p4) and (sp3)-(sp5) properties as in \S\ref{1.2}, but 
with thick ideals replaced by Serre ideals. In other words, these kinds of ideals can be defined entirely based on the notion of Serre ideals of abelian 2-categories, just like the more 
general prime ideals are defined in terms of thick ideals.

The set of Serre prime ideals of $\mc{T}$, denoted by $\mathrm{Serre}\mbox{-}\mathrm{Spec} (\mc{T})$, has an induced topology from 
$\Spec (\mc{T})$. This topology is shown to be intrinsically given in terms of Serre ideals of $\mc{T}$. If $\mc{C}$ is a strict abelian monoidal category, an alternative topology which more closely resembles the topology of Balmer in \cite{Balmer1} can also be put on $\mathrm{Serre}\mbox{-}\mathrm{Spec} (\mc{C})$. Under this topology, $\mathrm{Serre}\mbox{-}\mathrm{Spec} (\mc{C})$ is a ringed space.

Denote $\Zset_+ := \{ 0, 1, \ldots \}$. 
The Grothendieck ring $K_0(\mc{T})$ of an abelian 2-category $\mc{T}$ is a {\em{$\Zset_+$-ring}} in the terminology of \cite[Ch. 3]{Etingof1}, 
see Definitions \ref{dcate} and \ref{dZ+rings}. In \S \ref{Z+rings} we define the notions of {\em{Serre ideals}} and {\em{Serre prime}} 
({\em{semiprime}} and {\em{completely prime}}) {\em{ideals}} of a {\em{$\Zset_+$-ring $R$}}.
The set of Serre prime ideals of $R$, denoted by  $\mathrm{Serre}\mbox{-}\mathrm{Spec} (R)$, is equipped with a Zariski type topology.

It is proved in Section \ref{6} that, for an abelian 2-category $\mc{T}$ with the property that every 1-morphism of $\mc{T}$ has finite length, the 
functor $K_0$ induces bijections between the sets of Serre ideals, Serre prime (semiprime and completely prime ideals) of 
the abelian 2-category $\mc{T}$ and the $\Zset_+$-ring $K_0(\mc{T})$. Furthermore, the map
\[
K_0 : \mathrm{Serre}\mbox{-}\mathrm{Spec} (\mc{T}) \to \mathrm{Serre}\mbox{-}\mathrm{Spec} (K_0(\mc{T}))
\]
is shown to be a homeomorphism.

For a Serre prime (resp. semiprime, completely prime) ideal $\mc{I}$ of $\mc{T}$, the Serre quotient $\mc{T}/\mc{I}$ 
is a prime (resp. semiprime, domain) abelian 2-category. If every 1-morphism of $\mc{T}$ has finite length, then 
\[
K_0(\mc{T} / \mc{I}) \cong K_0(\mc{T}) / K_0(\mc{I}). 
\]

{\em{The point now is that if we have a categorification of a $\Zset_+$-ring $R$ via an abelian 2-category $\mc{T}$ (i.e., $K_0(\mc{T}) \cong R$)
and $I$ is a Serre ideal of $R$, then there is a unique Serre ideal $\mc{I}$ of $\mc{T}$ such that $K_0(\mc{I}) = I$. Furthermore,
the Serre quotient $\mc{T}/\mc{I}$ categorifies the $\Zset_+$-ring $R/I$, i.e., $K_0(\mc{T} / \mc{I}) \cong K_0(\mc{T}) / K_0 (\mc{I})$. If 
$I$ is a Serre prime (resp. semiprime, completely prime) ideal of the $\Zset_+$-ring $R$, then 
the Serre quotient $\mc{T}/\mc{I}$ is a prime (resp. semiprime, domain) abelian 2-category.}} We view this construction 
as a general way of constructing monoidal categorifications of $\Zset_+$-rings out of known ones by taking Serre quotients. 
This is illustrated in Section \ref{9} in the case of the {\em{quantized coordinate rings of open Richardson varieties}} for 
symmetric Kac--Moody algebras.

We expect that, in addition, Serre prime ideals of abelian 2-categories and $\Zset_+$-rings will provide a framework for finding intrinsic connections
between prime ideals of noncommutative algebras and totally positive parts of algebraic varieties. In the case of the algebras of quantum matrices, 
such a connection was previously found by exhibiting related explicit generating sets for prime ideals of the noncommutative algebras and 
minors defining totally positive cells \cite{GLL}.

{\em{Primitive ideals}} of abelian 2-categories $\mc{T}$ are introduced in Section \ref{7} as the annihilation ideals of simple exact 2-representations
in the setting of \cite{Mazorchuk2,Mazorchuk3}, where is proved that all such ideals are Serre prime ideals of $\mc{T}$. 
\subsection{Prime spectra of additive 2-categories}
One can develop analogous (but much simpler) theory of different forms of prime ideals of an {\em{additive 2-category}} $\mc{T}$, which is a 
2-category such that $\mc{T}(A_1, A_2)$ are 
additive categories for $A_1, A_2 \in \mc{T}$ and the compositions 
\[
\mc{T}(A_2, A_3) \times \mc{T}(A_1, A_2) \to \mc{T}(A_1, A_3)
\] 
are additive bifunctors for $A_1, A_2, A_3 \in \mc{T}$. 

This can be done by following exactly the same route as Sections \ref{3}-\ref{5} but based off the notion of a {\em{thick ideal
of an additive 2-category}}. Call a full subcategory of an additive category {\em{thick}} if its closed under direct sums, direct 
summands and isomorphisms. A {\em{thick ideal}} $\mc{I}$ of an additive 2-category $\mc{T}$ is a collection of subcategories $\mc{I}(A_1, A_2)$ of 
$\mc{T}(A_1, A_2)$ for all objects $A_1, A_2$ of $\mc{T}$ such that 
\begin{enumerate}
\item $\mc{I}(A_1, A_2)$ are thick subcategories of the additive categories $\mc{T}(A_1, A_2)$ and 
\item $\mc{T}_1 \circ \mc{I}_1 \subset \mc{I}_1$, $\mc{I}_1 \circ \mc{T}_1 \subseteq \mc{I}_1$. 
\end{enumerate}
Using the same conditions on containments with respect to thick ideals and 1-morphisms 
as in Sections \ref{3}-\ref{5}, one defines {\em{prime, semiprime, and completely prime ideals of 
additive 2-categories}} and proves analogs of the results in those sections (though in a simpler way than the abelian setting). 
One also analogously defines a Zariski topology on the set $\Spec(\mc{T})$ of prime ideals of $\mc{T}$ by using containments 
with respect to thick subcategories of $\mc{T}$.
There are no analogs of the Serre type ideals in this setting.

If an additive 2-category $\mc{T}$ has the property that each of its 1-morphisms has a unique decomposition as a finite set of indecomposables 
(e.g., if all additive categories $\mc{T}(A_1, A_2)$ are Krull--Schmidt), then the split Grothendieck ring $K_0^{sp}(\mc{T})$ is a $\Zset_+$-ring, 
see \reref{add}. Similarly to \S \ref{Zplus}, for such additive 2-categories $\mc{T}$, one shows that the map $K_0^{sp}(-)$ gives 
\begin{itemize}
\item a bijection between the sets of thick, prime, semiprime and completely prime ideals of $\mc{T}$ and the sets of Serre ideals, Serre prime, semiprime and completely 
prime ideals of the $\Zset_+$-ring $K_0^{sp}(\mc{T})$, and 
\item a homeomorphism $\Spec (\mc{T}) \to \mathrm{Serre}\mbox{-}\mathrm{Spec} (K_0^{sp}(\mc{T}))$.
\end{itemize}

Call a {\em{primitive ideal of an additive 2-category}} $\mc{T}$ to be the annihilation ideal of a simple 2-representation of $\mc{T}$ in the setting of 
\cite{Mazorchuk2,Mazorchuk3}. Similarly to Section \ref{7}, one shows that each such ideal is a prime ideal of $\mc{T}$.  

In a forthcoming publication we obtain analogs of the results in the paper for (noncommutative) prime spectra of triangulated 2-categories.
\subsection{Categorifications of Richardson varieties via prime Serre quotients} 
We finish with an important example of Serre completely prime ideals of abelian 2-categories that can be used to categorify 
the quantized coordinate rings of certain closures of open Richardson varieties. For a symmetrizable Kac--Moody group $G$, a pair of 
opposite Borel subgroups $B_\pm$ and Weyl group elements $u \leq w$, the corresponding open Richardson variety is defined as 
the intersections of opposite Schubert cells in the full flag variety of $G$,
\[
R_{u,w}:=(B_- u B_+)/B_+ \cap (B_+ w B_+)/B_+ \subset G/B_+.
\]
They have been used in a wide range of settings in representation theory, Schubert calculus, total positivity, Poisson geometry, and mathematical physics.
For symmetric Kac--Moody groups, Leclerc \cite{Leclerc} constructed a cluster algebra inside the coordinate ring of each Richardson variety of the same dimension. 
In the quantum situation, Lenagan and the second named author constructed large families of toric frames for all quantized coordinate rings of Richardson varieties 
that generate those rings \cite{LenYak}.

Recently, for each symmetrizable Kac--Moody algebra $\g$, 
Kashiwara, Kim, Oh, and Park \cite{KKOP} constructed a monoidal categorification of the quantization of a closure of $R_{u,w}$ in terms of 
a monoidal subcategory of the category of graded, finite dimensional representations of the 
Khovanov--Lauda--Rouquier (KLR) algebras associated to $\g$. Their construction uses Leclerc's interpretation of the coordinate ring of a closure of $R_{u,w}$ 
in terms of a double invariant subalgebra.

Denote by $\overline{R}_{u,w}$ the closure of $R_{u,w}$ in the Schubert cell $(B_+ w B_+)/B_+ \subset G/B_+$. We construct a monoidal 
categorification of the quantization $U^-_u[w]/I_u(w)$ of the coordinate ring of $\overline{R}_{u,w}$ used in the construction 
of toric frames in \cite{LenYak}. Here $U^-_q[w]$ are the quantum Schubert cell algebras \cite{DKP,Lusztig2} and $I_u(w)$ are the homogeneous 
completely prime ideals of these algebras that arose in the classification of their prime spectra in \cite{Y-plms}. This classification was based on the 
fundamental works of Anthony Joseph on the spectra of quantum groups \cite{Joseph-JA,Joseph-book} from the early 90s. It was proved 
in \cite{KKKO,Khovanov1,Rouquier0} that certain monoidal subcategories $\mc{C}_w$ of the categories of graded, finite dimensional modules 
of the KLR algebras associated to $\g$ categorify the dual integral form $U^-_\AA[w]\spcheck$ where $\AA:= \Zset[q^{\pm1}]$. 
We prove that for a symmetrizable Kac--Moody algebra $\g$, the ideals $I_w(u) \cap U^-_\AA[w]\spcheck$ have bases that are subsets 
of the upper global/canonical basis of $U^-_\AA[w]\spcheck$. From this we deduce that for symmetric $\g$,  $I_w(u) \cap U^-_\AA[w]\spcheck$ are Serre 
completely prime ideals of the $\Zset_+$-ring $U^-_\AA[w]\spcheck$. The bijection from \S \ref{1.3} implies that the monoidal category $\mc{C}_w$ 
has a Serre completely prime ideal $\mc{I}_u(w)$ such that $K_0(\mc{I}_u(w)) = I_w(u)$, and thus, the Serre quotient $\mc{C}_w / \mc{I}_u(w)$ categorifies 
$U^-_\AA[w]\spcheck / ( I_w(u) \cap U^-_\AA[w]\spcheck)$:
\[
K_0 (\mc{C}_w / \mc{I}_u(w)) \cong U^-_\AA[w]\spcheck / ( I_w(u) \cap U^-_\AA[w]\spcheck).
\]
It is an important problem to connect the categorification of Kashiwara, Kim, Oh, and Park \cite{KKOP} of open Richardson varieties (via subcategories of KLR modules) 
to ours (via Serre quotients of categories of KLR modules).
\medskip
\\
\noindent
{\bf Acknowledgements.} We are grateful to Masaki Kashiwara for communicating to us his proof of \thref{basis2} and to 
Ken Goodearl, Birge Huisgen-Zimmermann,
Osamu Iyama, Peter J{\o}rgensen, Bernhard Keller, Bernard Leclerc, Vanessa Miemietz and Michael Wemyss
for helpful discussions on thick subcategories of abelian and triangulated categories, 
canonical bases, and small 2-categories. We are thankful to the referee for the thorough reading of the paper and for making a number of valuable suggestions which helped us 
to improve the paper.

The research of K.V. was supported by a Board of Regents LSU fellowship and NSF grant DMS-1601862.
The research of M.Y. was supported by NSF grant DMS-1601862 and Bulgarian Science Fund grant DN 02-5.
\sectionnew{Abelian 2-Categories and categorification}
\label{two}
This section contains background material on (abelian) 2-categories and categorification of algebras.
\subsection{2-Categories}
\label{2-cat}
A category $\mc{T}$ is said to be {\em{enriched over a monoidal category}} $\mc{M}$ if the space of morphisms between any two objects of $\mc{T}$ is an object in 
$\mc{M}$ and $\mc{T}$ satisfies natural axioms which relate composition of morphisms in $\mc{T}$ and the identity morphisms of objects of $\mc{T}$
to the monoidal structure of $\mc{M}$. We refer the reader to \cite{Kelly1} for details. 

A {\em{2-category}} is a category enriched over the category of  1-categories. 
This means that for a 2-category $\mc{T}$, given two objects $A_1, A_2$ of it, the morphisms $\mc{T}(A_1, A_2)$ form a 1-category. 
The objects of these categories are denoted by the same symbol $\mc{T}(A_1, A_2)$ -- they are the 1-morphisms of $\mc{T}$. The morphisms 
of the categories $\mc{T}(A_1, A_2)$ are the 2-morphisms of $\mc{T}$. 
For a pair of 1-morphisms $f, g \in  \mc{T}(A_1, A_2)$, we will denote by $\mc{T}(f, g)$ the 2-morphisms between $f$ and $g$, i.e., 
the morphisms between the objects $f$ and $g$ in the category $\mc{T}(A_1, A_2).$

We have 2 types of compositions of 1- and 2-morphisms. We follow the notation of \cite{Lauda1}:
\begin{enumerate}
\item For a pair of objects $A_1, A_2$ of $\mc{T}$, the composition of morphisms in the category $\mc{T}(A_1, A_2)$ is called 
{\em{vertical composition}} of 2-morphisms of $\mc{T}$. In the {\em{globular representation}} of $\mc{T}$, such a
composition is given by the following diagram
\begin{center}
\begin{tikzpicture}
\node (a) at (0,0) {$A_2$};
\node (b) at (5,0) {$A_1.$};
\draw[<-] (a) to [bend left] node[scale=.7] (f) [above] {$h$} (b);
\draw[<-] (a) to node[scale=.7](g) [below] {$g$} (b);
\draw[<-] (a) to [bend right] node[scale=.7] (h) [below] {$f$} (b);
\draw[-{Implies},double distance=1.5pt,shorten >=2pt,shorten <=2pt] (g) to node[scale=.7] [right] {$\beta$} (f);
\draw[-{Implies},double distance=1.5pt,shorten >=2pt,shorten <=2pt] (h) to node[scale=.7] [right] {$\alpha$} (g);
\end{tikzpicture}
\end{center}
The vertical composition of the 2-morphisms $\alpha \in \mc{T}(f, g)$ and $\beta \in \mc{T}(g, h)$ will be denoted by 
$\beta \alpha \in \mc{T}(f, h)$, where $f, g, h$ are objects of $\mc{T}(A_1, A_2)$.

\item For each three objects $A_1, A_2, A_3$ of $\mc{T}$, we have a bifunctor of 1-categories
\begin{equation}
\label{composition}
\mc{T}(A_2, A_3) \times \mc{T}(A_1, A_2) \to \mc{T}(A_1, A_3).
\end{equation}
The resulting composition of 1- and 2-morphisms of $\mc{T}$ is called  {\em{horizontal composition}}. 
In the globular representation of $\mc{T}$, these compositions are given by the diagram
\begin{center}
\begin{tikzpicture}
\node (a) at (0,0) {$A_3$};
\node (b) at (4,0) {$A_2$};
\node (c) at (8,0) {$A_1.$};
\draw[<-] (a) to [bend left] node[scale=.7] (f) [above] {$f_2$} (b);
\draw[<-] (a) to [bend right] node[scale=.7] (g) [below] {$g_2$} (b);
\draw[-{Implies},double distance=1.5pt,shorten >=2pt,shorten <=2pt] (g) to node[scale=.7] [right] {$\alpha_2$} (f);

\draw[<-] (b) to [bend left] node[scale=.7] (f) [above] {$f_1$} (c);
\draw[<-] (b) to [bend right] node[scale=.7] (g) [below] {$g_1$} (c);
\draw[-{Implies},double distance=1.5pt,shorten >=2pt,shorten <=2pt] (g) to node[scale=.7] [right] {$\alpha_1$} (f);
\end{tikzpicture}
\end{center}
In this notation, the horizontal composition of 2-morphisms will be denoted by $\alpha_2 *  \alpha_1$. The 
horizontal composition of 1-morphisms will be denoted by $f_2 f_1$. 
\end{enumerate}
A 2-category $\mc{T}$ has identity 1-morphisms $1_A \in \mc{T}(A,A)$ (for its objects $A \in \mc{T}$). The compositions and
identity morphisms satisfy natural associativity and identity axioms \cite{Lauda1,Mazorchuk1}, 
which are equivalent to the definition of 2-categories in the language of enriched categories.

{\em{2-categories}} are generalizations of {\em{monoidal categories}}, in the sense that a strict monoidal category is the same thing as a 2-category with one object:

To a strict monoidal category $\mc{M}$, one associates a 2-category $\mc{T}$ with one object $A$ by taking $\mc{T}(A, A):= \mc{M}$. The tensor product in 
$\mc{M}$ is used to define composition of 1-morphisms of $\mc{T}$. For $f,g \in \mc{M}= \mc{T}(A,A)$, one defines the 2-morphisms $\mc{T}(f,g):= \mc{M}(f,g)$. 
All 2-categories with 1 object arise in this way.

Recall that a 1-category $\mc{C}$ is called {\em{small}} 
if its objects form a set and $\mc{C}(A_1, A_2)$ is a set for all pairs of objects $A_1, A_2 \in \mc{C}$.
{\em{Throughout the paper we work with small 2-categories $\mc{T}$, which are 2-categories satisfying the conditions that 
the objects of $\mc{T}$ form a set and $\mc{T}(A_1, A_2)$ is a small 1-category for all pairs of objects $A_1, A_2$ of $\mc{T}$.}}

The set of objects of such a 2-category $\mc{T}$ will be denoted by the same symbol $\mc{T}$. The set of 
1-morphisms of $\mc{T}$ will be denoted by $\mc{T}_1$.
\subsection{Abelian 2-categories and categorification}
\label{ab-2-cat}
\bde{ab2cat} We will say that a 2-category $\mc{T}$ is an {\em{abelian 2-category}} if $\mc{T}(A_1, A_2)$ are 
abelian categories for all $A_1, A_2 \in \mc{T}$ and the compositions 
\[
\mc{T}(A_2, A_3) \times \mc{T}(A_1, A_2) \to \mc{T}(A_1, A_3)
\] 
are exact bifunctors for all $A_1, A_2, A_3 \in \mc{T}$.

More generally, for a ring $\kk$, we will say that $\mc{T}$ is a {\em{$\kk$-linear abelian 2-category}}
if $\mc{T}(A_1, A_2)$ are $\kk$-linear abelian categories for $A_1, A_2 \in \mc{T}$.
\ede

A {\em{multiring category}} in the terminology of \cite[Definition 4.2.3]{Etingof1} is precisely a 
{\em{$\kk$-linear abelian 2-category with one object}}. 

\bre{enrich-ab-2}
Let $\kk$ be a field. Recall that a $\kk$-linear abelian category $\mc{C}$ is called {\em{locally finite}} if it is $\Hom$-finite 
(i.e.,  $\dim_{\kk} \mc{C}(A_1, A_2) < \infty$ for all $A_1, A_2 \in \mc{C})$ and each object of $\mc{C}$ has finite length; 
we refer the reader to \cite[\S 1.8]{Etingof1} for details. 
Let $\mc{LF}Ab_{ex}$ be the monoidal category of locally finite abelian categories equipped with the Deligne tensor product
(\cite{Deligne} and \cite[\S 1.11]{Etingof1}) and morphisms given by exact functors. 

In this terminology, 
a $\kk$-linear abelian 2-category $\mc{T}$ with the property that the 1-categories $\mc{T}(A_1, A_2)$ are locally finite 
for all $A_1, A_2 \in \mc{T}$ is the same thing as a category which is enriched over the monoidal category $\mc{LF}Ab_{ex}$.
This is easy to verify, the only key step being the universality property of the Deligne tensor product with respect to exact functors
\cite[Proposition 1.11.2(v)]{Etingof1}.
\ere

We will denote by $K_0(\mc{C})$ the Grothendieck group of an abelian category $\mc{C}$. To each abelian 2-category $\mc{T}$ one associates the 
pre-additive category $K_0(\mc{T})$ whose objects are the objects of $\mc{T}$ and morphisms are
\[
K_0(\mc{T})(A_1, A_2) := K_0( \mc{T}(A_1, A_2)) \quad \mbox{for} \; \; A_1, A_2 \in \mc{T}. 
\]
Given a pre-additive category $\mc{F}$, one says that the 2-category $\mc{T}$ categorifies $\mc{F}$ 
if $K_0(\mc{T}) \cong \mc{F}$ as pre-additive categories.

To a pre-additive category $\mc{F}$, one associates a ring with elements
\[
\oplus_{A_1, A_2 \in \mc{F} } \mc{F}(A_1, A_2).
\]
The product in the ring is the composition of morphisms when it makes sense and 0 otherwise. 
In particular, the identity morphisms $1_A$ are idempotents of the ring for all objects $A \in \mc{F}$. By abuse of notation, this ring 
is denoted by the same symbol $\mc{F}$ as the original category.

\bde{cate} For an abelian 2-category $\mc{T}$, the ring $K_0(\mc{T})$ is called the Grothendieck ring of $\mc{T}$.
We say that $\mc{T}$ categorifies an $S$-algebra $R$, for a commutative ring $S$, if $K_0(\mc{T}) \otimes_{\mathbb{Z}} S \cong R$.
\ede

Often, it is not sufficient to consider multiring categories (abelian monoidal categories) to obtain categorifications of algebras, and one needs 
the more general setting of 2-categories.

\bre{add} An additive 2-category is a 2-category $\mc{T}$ such that $\mc{T}(A_1, A_2)$ are 
additive categories for all $A_1, A_2 \in \mc{T}$ and the compositions 
\[
\mc{T}(A_2, A_3) \times \mc{T}(A_1, A_2) \to \mc{T}(A_1, A_3)
\] 
are additive bifunctors for all $A_1, A_2, A_3 \in \mc{T}$. For such a category $\mc{T}$, one defines the 
pre-additive category $K_0^{sp}(\mc{T})$ whose objects are the objects of $\mc{T}$ and morphisms are the split Grothendieck groups 
\[
K_0^{sp} (\mc{T} (A_1, A_2))
\] 
of the additive categories $\mc{T}(A_1, A_2)$ for $A_1, A_2 \in \mc{T}$. 

We say that an additive 2-category $\mc{T}$ categorifies an $S$-algebra $R$ if $K_0^{sp}(\mc{T}) \otimes_{\mathbb{Z}} S \cong R$.
\ere
\sectionnew{The prime spectrum}
\label{3}
In this section we define the prime spectrum of an abelian 2-category and a Zariski type topology on it. 
We prove two equivalent characterizations of prime ideals, extending theorems from classical ring theory.
We also prove that maximal elements of the sets of ideals not intersecting multiplicative sets 
of 1-morphisms of 2-categories are prime ideals. 
\subsection{Thick ideals of abelian 2-categories}
\label{thick-2-cat}
\bde{weak-sub} A {\em{weak subcategory}} $\mc{I}$ of a 2-category $\mc{T}$ is
\begin{enumerate}
\item a subcollection $\mc{I}$ of objects of $\mc{T}$ and 
\item a collection of subcategories $\mc{I}(A_1, A_2)$ of $\mc{T}(A_1, A_2)$ for $A_1, A_2 \in \mc{I}$,
\end{enumerate}
such that the composition bifunctors \eqref{composition} restrict to bifunctors
\[
\mc{I}(A_2, A_3) \times \mc{I}(A_1, A_2) \to \mc{I}(A_1, A_3)
\] 
for $A_1, A_2, A_3 \in \mc{I}$.
\ede
A weak subcategory $\mc{I}$ of a 2-category $\mc{T}$ is not necessarily a 2-category on its own because it might not contain the identity 
morphisms $1_A$ for its objects $A \in \mc{I}$. Apart from this, a weak subcategory of a 2-category satisfies the other axioms for 2-categories.
The relationship of a weak subcategory to a 2-category is the same as the relationship of a subring to a unital ring $R$. In the latter case, the 
subring does not need to contain the unit of $R$.

\bde{thick-id} (1) A {\em{thick subcategory}} of an abelian category is a nonempty full subcategory which is closed under taking kernels, 
cokernels, and extensions.

(2) A {\em{thick weak subcategory}} of an abelian 2-category $\mc{T}$ is a weak subcategory $\mc{I}$ of $\mc{T}$ having the same set of objects and
such that for any pair of objects $A_1, A_2 \in \mc{T}$, $\mc{I}(A_1, A_2)$ is a thick subcategory of the abelian category $\mc{T}(A_1,A_2)$. 

(3) A {\em{thick ideal}} of an abelian 2-category $\mc{T}$ is a thick weak subcategory $\mc{I}$ of $\mc{T}$ such that, 
for all 1-morphisms $f$ in $\mc{T}$ and $g$ in $\mc{I}$, the compositions $fg, gf$ are 1-morphisms of $\mc{I}$ 
whenever they are defined.
\ede

Sometimes the term {\em{wide subcategory}} is used instead of {\em{thick}}, see for instance \cite{IT}.

Every thick subcategory of an abelian category is closed under isomorphisms and taking direct summands of its objects 
(because one can take the kernels of idempotent endomorphisms of its objects).

For a thick weak subcategory $\mc{I}$ of an abelian 2-category $\mc{T}$, $\mc{I}(A_1,A_2)$ is an abelian category for 
every pair of objects $A_1, A_2 \in \mc{I}$ with respect to the same kernels and cokernels as the ambient abelian category 
$\mc{T}(A_1,A_2)$.

In part (3), the compositions of 1-morphisms that are used are the horizontal compositions discussed in \S \ref{2-cat}.
More explicitly, a thick subcategory $\mc{I}$ of $\mc{T}$ is a thick ideal if for all $f_1 \in \mc{T}(A_1,A_2)$,
$g_2 \in \mc{I}(A_2,A_3)$ and $f_3 \in \mc{T}(A_3,A_4)$, we have
\[
g_2 f_1 \in  \mc{I}(A_1,A_3) \quad \mbox{and} \quad f_3 g_2 \in  \mc{I}(A_2,A_4),
\]
where $A_1, A_2, A_3, A_4 \in \mc{T}$.

\bre{incl-thick} Let $\mc{I}$ and $\mc{J}$ be a pair of thick weak subcategories of $\mc{T}$. Then
\[
\mc{I} \subseteq \mc{J} \quad \mbox{if and only if} \quad \mc{I}_1 \subseteq \mc{J}_1.
\]
In particular,
\[
\mc{I} = \mc{J} \quad \mbox{if and only if} \quad \mc{I}_1 = \mc{J}_1.
\]
\ere

\bex{0} There exists a unique thick ideal of every 2-category $\mc{T}$ whose set of 1-morphisms consists of the 0 objects of the abelian 
categories $\mc{T}(A_1, A_2)$ for all $A_1, A_2 \in \mc{T}$. This thick ideal will be denoted by $0_{\mc{T}}$. Every other 
thick ideal of $\mc{T}$ contains $0_{\mc{T}}$. 
\eex

For two subsets $X, Y \subseteq \mc{T}_1$, denote by 
\begin{align*}
\mbox{$X \circ Y$} \; &\mbox{the set of isomorphism classes of 1-morphims of $\mc{T}$ having representatives of the}
\\
&\mbox{form $f g$ for $f$ and $g$ representing classes in $X$ and $Y$ such that $fg$ is defined}.
\end{align*}
In general, $X \circ Y$ can be empty. For $f, g \in \mc{T}_1$ the composition $f \circ g$ is either empty or consists of one element.

In this notation, a thick weak subcategory $\mc{I}$ of $\mc{T}$ is a thick ideal if and only if
\[
\mc{T}_1 \circ \mc{I}_1 \subseteq \mc{I}_1 \quad \mbox{and} 
\quad \mc{I}_1 \circ \mc{T}_1 \subseteq \mc{I}_1.
\]

\bde{1-sided} A {\em{thick left (respectively right) ideal}} of an abelian 2-category $\mc{T}$ is a thick weak subcategory $\mc{I}$ of $\mc{T}$ such that
\[
\mc{T}_1 \circ \mc{I}_1 \subseteq \mc{I}_1 \quad \mbox{(respectively} 
\; \; 
\mc{I}_1 \circ \mc{T}_1 \subseteq \mc{I}_1).
\]
\ede

\bre{isom} Let $A_1, A_2, B_1, B_2$ be four objects of an abelian 2-category $\mc{T}$ such that 
\[
A_i \cong B_i \quad \mbox{for} \quad i=1,2. 
\]
Then for every thick ideal $\mc{I}$ of $\mc{T}$, we have (noncanonical) isomorphisms of abelian categories
\begin{equation}
\label{iso-IAB}
\mc{I} (A_1, A_2) \cong \mc{I} (B_1, B_2).
\end{equation}
Indeed, let $f_i \in \mc{T}(A_i, B_i)$ and $g_i \in \mc{T}(B_i, A_i)$ be such that 
\[
f_i g_i \cong 1_{B_i} \quad \mbox{and} \quad \quad g_i f_i \cong 1_{A_i}
\]
for $i=1,2$ (where the isomorphisms are in the categories $\mc{T}(B_i,B_i)$ and $\mc{T}(A_i,A_i)$).
The functor giving the equivalence \eqref{iso-IAB} is defined by 
\[
h \mapsto f_2 h g_1
\]
on the level of objects $h \in \mc{I}(A_1, A_1)$ and 
\[
\alpha \mapsto 1_{f_2} * \alpha * 1_{g_1}
\]
on the level of morphisms.
\ere
\subsection{Prime ideals of abelian 2-categories}
\label{prime-2-cat}
A thick ideal $\mc{I}$ of $\mc{T}$ will be called {\em{proper}} if $\mc{I} \neq \mc{T}$; by \reref{incl-thick} this is the same as 
$\mc{I}_1 \subsetneq \mc{T}_1$.

\bde{prime-ideals}
We call $\mc{P}$ a {\em{prime ideal}} of $\mc{T}$ if $\mc{P}$ is a proper thick ideal of $\mc{T}$ with the property that 
for every pair of thick ideals $\mc{I}$ and $\mc{J}$ of $\mc{T}$, 
\[
\mc{I}_1 \circ \mc{J}_1 \subseteq \mc{P}_1 \quad \Rightarrow \quad \mc{I} \subseteq \mc{P} \; \;  \mbox{or} \; \;  \mc{J} \subseteq \mc{P}. 
\]
The set of all prime ideals $\mc{P}$ of an abelian 2-category $\mc{T}$ will be called the {\em{prime spectrum}} of $\mc{T}$ and will be denoted by $\Spec(\mc{T})$. 
\ede
By \reref{incl-thick}, the property on the right side of the implication can be replaced with $\mc{I}_1 \subseteq \mc{P}_1$ or $\mc{J}_1 \subseteq \mc{P}_1$.
\subsection{Two equivalent characterizations of prime ideals}
\label{equiv-prim-id}
The following lemma is straightforward.
\ble{int-thick}
The intersection of any family of thick ideals is a thick ideal. 
\ele

If $\mc{M}$ is a collection of 1-morphisms of $\mc{T}$ (i.e., $\mc{M} \subseteq \mc{T}_1$), let $\langle \mc{M} \rangle$ denote the smallest thick ideal of $\mc{T}$ containing $\mc{M}$,
which exists by the previous lemma.

\ble{MTN} For every two collections $\mc{M}, \mc{N} \subseteq \mc{T}_1$ of 1-morphisms of an abelian 2-category $\mc{T}$, we have
\begin{equation}
\label{MTN}
\langle \mc{M} \rangle_1 \circ \langle \mc{N} \rangle_1 \subseteq \langle \mc{M} \circ \mc{T}_1 \circ \mc{N} \rangle_1.
\end{equation}
\ele
\begin{proof}
We will first show that 
\begin{equation}
\label{first-emb}
\langle \mc{M} \rangle_1 \circ \mc{N} \subseteq \langle \mc{M} \circ \mc{T}_1 \circ \mc{N} \rangle_1.
\end{equation}
The 1-morphisms of $\langle \mc{M} \rangle$ are obtained from the elements of $\mc{M}$ by successive taking of kernels and
cokernels (of 2-morphisms between these elements), and extensions (between these elements), 
as well as compositions on 
the left and the right by elements in $\mc{T}_1$. 
We need to show that those operations, composed on the right 
with the elements of $\mc{N}$, yield elements of the right hand side.

(1) Suppose that $\alpha: f \to g$ is a 2-morphism for $f, g \in \mc{T}_1$ with the property that 
\[
f n, g n \in \langle \mc{M} \circ \mc{T} \circ \mc{N} \rangle_1 \quad \mbox{for all} \quad n \in \mc{T}_1 \circ \mc{N}.
\]
Note that, for example, every 1-morphism in $\mc{M}$ has this property. 
Let $\kappa :  k \to f$ be the kernel of $\alpha$. Since exact functors preserve kernels,
$\kappa * \id_n: k n \to f n$ is the kernel of $\alpha * \id_n: f n \to g n$. The thickness property of 
$\langle \mc{M} \circ \mc{T}_1 \circ \mc{N} \rangle$ implies that $k n \in \langle \mc{M} \circ \mc{T}_1 \circ \mc{N} \rangle_1$ for all
$n \in \mc{T}_1 \circ \mc{N}$. 

Symmetrically, one shows that if $\gamma : g  \to c$ is the cokernel of $\alpha$, 
then $c n \in \langle \mc{M} \circ \mc{T}_1 \circ \mc{N} \rangle_1$ for all $n \in \mc{T}_1 \circ \mc{N}$.

(2) Next, assume that 
\[
0 \to f \to g \to h \to 0
\]
is an exact sequence in one of the abelian categories $\mc{T}(A_1, A_2)$, where $f, h$ have the property that 
$f n, h n \in \langle \mc{M} \circ \mc{T}_1 \circ \mc{N} \rangle_1$ for all $n \in \mc{T}_1 \circ \mc{N}$.
Since horizontal composition in $\mc{T}$ is exact, for any $n \in \mc{T}_1 \circ \mc{N}$,
we get a short exact sequence 
\[
0 \to f  n \to g n \to h n \to 0.
\]
Since the first and last terms are in $\langle \mc{M} \circ \mc{T}_1 \circ \mc{N} \rangle_1$, so is the middle term. 

Combining (1)--(2) and the fact that $\langle \mc{M} \circ \mc{T}_1 \circ \mc{N} \rangle_1$ is stable under left compositions with elements of $\mc{T}_1$, 
yields \eqref{first-emb}. Analogously, we derive \eqref{MTN} from \eqref{first-emb} by using $\langle \mc{M} \rangle_1$ in place of $\mc{M}$. 
\end{proof}

\bth{prime1} A proper thick ideal $\mc{P}$ of an abelian 2-category $\mc{T}$ is prime if and only if for all $m, n \in \mc{T}_1$, 
$m \circ \mc{T}_1 \circ n \subseteq \mc{P}_1$ implies that either $m \in \mc{P}_1$ or $n \in \mc{P}_1$.
\eth
\begin{proof} Suppose $\mc{P}$ is a prime ideal of $\mc{T}$, and that $m \circ \mc{T} \circ n \subseteq \mc{P}$ for some $m, n \in \mc{T}_1$.
Then by the previous lemma,
\[
\langle m \rangle_1 \circ \langle n \rangle_1 \subseteq \langle m \circ \mc{T}_1 \circ n \rangle_1 \subseteq \mc{P}_1,
\]
and so by primeness of $\mc{P}$, $\langle m \rangle \subseteq \mc{P}_1$ or $\langle n \rangle \subseteq \mc{P}_1$.
Therefore, $m$ or $n$ is in $\mc{P}_1$.

Now suppose $\mc{P}$ is a proper thick ideal of $\mc{T}$ with the property that for all $m, n \in \mc{T}_1$, 
$m \circ \mc{T}_1 \circ n \subseteq \mc{P}_1$ implies that either $m \in \mc{P}_1$ or $n \in \mc{P}_1$.
Let $\mc{I}$ and $\mc{J}$ be a pair of thick ideals of $\mc{T}$ such that 
\[
\mc{I}_1 \circ \mc{J}_1 \subseteq \mc {P}_1 \quad \mbox{and} \quad \mc{J}_1 \not \subseteq \mc{P}_1.
\]
Then there is some $j \in \mc{J}_1$ with $j \not \in \mc{P}_1$.
However, $i \circ \mc{T}_1 \circ j \subseteq \mc{P}_1$ for any $i \in \mc{I}_1$ 
since $i \circ \mc{T}_1 \subseteq \mc{I}_1, j \in \mc{J}_1$, and $\mc{I}_1 \circ \mc{J}_1 \subseteq \mc{P}_1$. 
The assumed property of $\mc{P}$ implies that $\mc{I}_1 \subseteq \mc{P}_1$. 
Therefore $\mc{I} \subseteq \mc{P}$ by \reref{incl-thick}.
\end{proof}

It is easy to see that \thref{prime1} implies the following:
\bpr{prime1b} A proper thick ideal $\mc{P}$ of an abelian 2-category $\mc{T}$ is prime if and only if  for 
every pair of right thick ideals $\mc{I}$ and $\mc{J}$ of $\mc{T}$ 
\[
\mc{I}_1 \circ \mc{J}_1 \subseteq \mc{P}_1 \quad \Rightarrow \quad \mc{I} \subseteq \mc{P} \; \;  \mbox{or} \; \;  \mc{J} \subseteq \mc{P}. 
\]
\epr
A similar characterization holds using left thick ideals.

\bth{prime2} A proper thick ideal $\mc{P}$ of an abelian 2-category $\mc{T}$ is prime if and only if for all thick ideals $\mc{I}, \mc{J}$ of $\mc{T}$ 
properly containing $\mc{P}$, we have that $\mc{I}_1 \circ \mc{J}_1 \not \subseteq \mc{P}_1$. 
\eth
\begin{proof} The implication $\Rightarrow$ is clear. Suppose $\mc{P}$ is a proper thick ideal which is not prime. Then there exist some thick ideals $\mc{I}$ and $\mc{J}$ of $\mc{T}$ 
with $\mc{I}_1 \circ \mc{J}_1 \subseteq \mc{P}_1$ and $\mc{I}, \mc{J} \not \subseteq \mc{P}$. Set
\[
\mc{M} := \mc{P}_1 \cup \mc{I}_1 \quad \mbox{and} \quad  \mc{N} := \mc{P}_1 \cup \mc{J}_1.
\]  
By \reref{incl-thick}, $\mc{P}_1$ is properly contained in both $\langle \mc{M} \rangle_1$ and $\langle \mc{N} \rangle_1$. \leref{MTN} implies that 
\begin{equation}
\langle \mc{M} \rangle \circ \langle \mc{N} \rangle \subseteq \langle \mc{M} \circ \mc{T}_1 \circ \mc{N} \rangle.
\label{bylem}
\end{equation}
Observe that 
\begin{equation}
\mc{M} \circ \mc{T}_1 \circ \mc{N} \subseteq \mc{P}_1, 
\label{itj}
\end{equation}
by the following. Consider the composition $i t  j$ for some $i \in \mc{M}, t \in \mc{T}_1, j \in \mc{N}$. 
So, $i \in \mc{I}$ or $\mc{P}_1$; likewise, $j \in \mc{I}$ or $\mc{P}_1$. If at least one of the two 1-morphism $i,j$ is in $\mc{P}_1$, we have $i t  j \in \mc{P}_1$
since $\mc{P}$ is a thick ideal; if $i \in \mc{I}_1$ and $j \in \mc{J}_1$ then $i \circ t \in \mc{I}_1$, so $i  t  j \in \mc{I}_1 \circ \mc{J}_1 \subseteq \mc{P}_1$ by assumption. 

Therefore  $\langle \mc{M} \rangle$ and $\langle \mc{N} \rangle$ are thick ideals properly containing $\mc{P}$ and 
$\langle \mc{M} \rangle_1 \circ \langle \mc{N} \rangle_1 \subseteq \mc{P}_1$
(the last inclusion follows from \eqref{bylem}--\eqref{itj} and the minimality of the thick ideal $\langle - \rangle$).
Hence, $\mc{P}$ does not have the stated property, 
which completes the proof of the theorem.
\end{proof}
\subsection{Relation to maximal ideals}
\label{max-id}
\bde{mult} A nonempty set  $\mc{M} \subseteq \mc{T}_1$ will be called {\em{multiplicative}} if $\mc{M}$ is a subset of non-zero equivalence classes of objects 
of $\mc{T}(A,A)$ for some object $A$ of $\mc{T}$ and $\mc{M} \circ \mc{M} \subseteq \mc{M}$. 
\ede

The condition that $\mc{M} \subseteq \mc{T}(A,A)$ means that all 1-morphism in $\mc{M}$ are composable. Let us explain the motivation for this condition. 
Let $R$ be a ring and $\{e_s\}$ be a collection of orthogonal idempotents. If $M$ is a multiplicative subset such that 
\[
M \subseteq \bigcup_{s,t} e_s R e_t
\]
then $M \subseteq e_s R e_s$ for some $s$, because otherwise $M$ will contain the 0 element of $R$.

\bth{mult-max} Assume that $\mc{M}$ is a multiplicative subset of $\mc{T}_1$ for an abelian 2-category $\mc{T}$ and that 
$\mc{I}$ is a proper thick ideal of $\mc{T}$ such that $\mc{I}_1 \cap \mc{M} = \varnothing$.

Let $\mc{P}$ be a maximal element of the 
collection of thick ideals of $\mc{T}$ containing $\mc{I}$ and intersecting $\mc{M}$ trivially, equipped with the inclusion relation, i.e., 
$\mc{P}$ is a maximal element of the set 
\[
X(\mc{M}, \mc{I}) := \{ \mc{K} \; \mbox{a thick ideal of} \; \; \mc{T} \mid \mc{K} \supseteq \mc{I} \; \; \mbox{and} \; \; \mc{K}_1 \cap \mc{M} = \varnothing \}.
\]
Then $\mc{P}$ is prime.
\eth

\begin{proof}
Fix such an ideal $\mc{P}$. 
Suppose $\mc{Q}$ and $\mc{R}$ are thick ideals properly containing $\mc{P}$. By \thref{prime2}, it is enough to show that $\mc{Q} \circ \mc{R} \not \subseteq \mc{P}$.
Since
\[
\mc{I} \subseteq \mc{P} \subseteq \mc{Q} \quad \mbox{and} \quad
\mc{I} \subseteq \mc{P} \subseteq \mc{R},
\]
both $\mc{Q}_1$ and $\mc{R}_1$ must intersect nontrivially with $\mc{M}$, by the maximality assumption on $\mc{P}$. 
Let $q \in \mc{Q}_1 \cap \mc{M}$ and $r \in \mc{R}_1 \cap \mc{M}$.
If $\mc{Q} \circ \mc{R} \subseteq \mc{P}$, then we would obtain that $q r \in \mc{P}$, because by the definition of multiplicative subset of $\mc{T}_1$, 
each two elements of $\mc{M}$ are composable. However, since $q  r \in \mc{M}$, this contradictions with the assumption that  $\mc{P}_1 \cap \mc{M} = \varnothing$. 
\end{proof}

\bre{max-M} The set $X(\mc{M}, \mc{I})$ from \thref{mult-max} is nonempty because $\mc{I} \in X(\mc{M}, \mc{I})$. 
The union of an ascending chain of thick ideals in the set $X(\mc{M}, \mc{I})$ is obviously a thick ideal of $\mc{T}$. 
By Zorn's lemma, the set $X(\mc{M}, \mc{I})$ from \thref{mult-max} always contains at least one maximal element.
\ere
\bco{nonempty} (1) For each proper thick ideal $\mc{I}$ of an abelian 2-category $\mc{T}$, there 
exists a prime ideal $\mc{P}$ of $\mc{T}$ that contains $\mc{I}$. 
 
(2) Let $\mc{M}$ be a multiplicative set of an abelian 2-category $\mc{T}$. Every maximal element of the 
set of thick ideals $\mc{K}$ of $\mc{T}$ such that $\mc{K}_1 \cap \mc{M} = \varnothing$ 
is a prime ideal. The set of such thick ideals contains at least one maximal element.
\eco
\begin{proof} (1) Since the thick ideal $\mc{I}$ is proper, there exists an object $A \in \mc{T}$ such that $1_A \notin \mc{I}_1$. Indeed, 
otherwise 
\[
\mc{T}(B,A) = \mc{T}(B,A) \circ 1_A = \mc{I}(B,A)
\]
for all objects $A,B \in \mc{T}$. The statement of part (1) follows from \thref{mult-max} applied to $\mc{M} := \{ 1_A \}$ 
for an object $A \in \mc{T}$ such that $1_A \notin \mc{I}$. 

(2) For each multiplicative subset $\mc{M}$ of an abelian 2-category $\mc{T}$, the thick ideal $0_{\mc{T}}$ 
from \exref{0} intersects $\mc{M}$ trivially. This part follows from \thref{mult-max} applied to the thick ideal
$\mc{I} := 0_{\mc{T}}$.
\end{proof}

The second part of the corollary, applied to the multiplicative subset $\mc{M}:= \{1_A\}$ for an object $A \in \mc{T}$, 
implies the following: 
 \bco{nonempS} 
The prime spectrum of every abelian 2-category $\mc{T}$ is nonempty. 
\eco
\bde{prime-simp} (1) An abelian 2-category $\mc{T}$ will be called {\em{prime}} if $0_{\mc{T}}$ is a 
prime ideal of $\mc{T}$. 

(2) An abelian 2-category $\mc{T}$ will be called {\em{simple}} if the only proper thick ideal of $\mc{T}$ is $0_{\mc{T}}$. 
\ede
\coref{nonempS} implies that every simple abelian 2-category $\mc{T}$ is prime.
\subsection{The Zariski topology}
\label{Zariski}
\bde{Za} Define the family of closed sets $V(\mc{I}):= \{ \mc{P} \in \Spec( \mc{T}) \mid \mc{P} \supseteq \mc{I} \}$ 
of $\Spec (\mc{T})$ for all thick ideals $\mc{I}$. 
\ede
\bre{different}
This topology is different from the one considered by Balmer \cite{Balmer1}. The main reason for which we 
consider it is to ensure good behavior under the $K_0$ map, see \thref{fin-ab}(4).
\ere
\ble{Zto} For each abelian 2-category $\mc{T}$, 
the above family of closed sets turns $\Spec ( \mc{T})$ into a topological space.  
The corresponding topology will be called the Zariski topology of $\Spec ( \mc{T})$.
\ele
It is easy to verify for that for every pair of  thick ideals $\mc{I}, \mc{J}$ of $\mc{T}$  and for every 
(possibly infinite) collection $\{\mc{I}_s\}$ of thick ideals of $\mc{T}$, similarly to the classical situation, we have
\begin{align*}
V(\mc{I}) \cup V ( \mc{J})&=V( \langle \mc{I}_1 \circ \mc{J}_1 \rangle ) \quad \mbox{and}
\\
\bigcap_i V(\mc{I}_s) &=V\left ( \left \langle \bigcup_i (\mc{I}_s)_1 \right \rangle \right ). 
\end{align*}
Finally, we also have  $V(\mc{T})= \varnothing$ and $V( 0_{\mc{T}} )= \Spec (\mc{T}).$ 
\subsection{An example}
\label{matr-ex}
Let $\Gamma$ be a nonempty set and $\kk$ be an arbitrary field. Let ${}_{\kk}\mathrm{Vect}_{\kk}$ be the category of finite dimensional $\kk$-vector spaces considered 
as $(\kk, \kk)$-bimodules. Let $\{ \kk_a \mid a \in \Gamma \}$ be a collection of fields isomorphic to $\kk$ and indexed by $\Gamma$. 

There is a unique $\kk$-linear abelian 2-category ${\mathscr{M}}_{\Gamma}(\kk)$ whose set of objects is $\Gamma$ and such that 
\[
{\mathscr{M}}_{\Gamma}(\kk) (a, b):= {}_{\kk_b} \mathrm{Vect}_{\kk_a} \quad \mbox{for} \quad a, b \in \Gamma.
\]
Its composition bifunctors are given by
\[
-  \otimes_{\kk_b} - : {\mathscr{M}}_{\Gamma}(\kk) (b, c) \times {\mathscr{M}}_{\Gamma}(\kk) (a, b) \to {\mathscr{M}}_{\Gamma}(\kk) (a, c).
\]
Its Grothendieck group is $K_0({\mathscr{M}}) \cong M_\Gamma(\Zset)$ -- the ring of square matrices with finitely many nonzero 
integer entries whose rows and columns are indexed by $\Gamma$. In the terminology of \S \ref{ab-2-cat}, 
${\mathscr{M}}_{\Gamma}(\kk)$ is a categorification of the matrix ring $M_\Gamma(\kk')$ for any field $\kk'$. 

Analogously to the classical situation, we show: 
\ble{mat-simp} The abelian 2-categories ${\mathscr{M}}_{\Gamma}(\kk)$
are simple (and thus prime).
\ele
\begin{proof} Let $\mc{I}$ be a thick ideal of ${\mathscr{M}}_{\Gamma}(\kk)$ that properly contains 
the 0-ideal $0_{{\mathscr{M}}_{\Gamma}(\kk)}$. Then for some $a, b \in \Gamma$, 
\[
\mc{I}(a,b) \neq 0.
\]
Since $\mc{I}$ is thick, $\mc{I}(a,b)$ is a nonzero subcategory of ${}_{\kk_b} \mathrm{Vect}_{\kk_a}$ that is closed under taking direct summands. 
Hence $\mc{I}(a,b)$ contains the 1-dimensional vector space in ${}_{\kk_b} \mathrm{Vect}_{\kk_a}$, and so, 
\[
\mc{I}(a,b) = {}_{\kk_b} \mathrm{Vect}_{\kk_a} = {\mathscr{M}}_{\Gamma}(\kk) (a, b). 
\]
Since all objects in ${\mathscr{M}}_{\Gamma}(\kk)$ are isomorphic to each other, \reref{isom} implies that 
\[
\mc{I}(a',b') = {\mathscr{M}}_{\Gamma}(\kk) (a', b')
\]
for all $a', b' \in \Gamma$. Thus $\mc{I} = {\mathscr{M}}_{\Gamma}(\kk)$, which completes the proof.
\end{proof}
\sectionnew{Minimal primes in noetherian abelian 2-categories}
\label{4}
In this section we define noetherian abelian 2-categories $\mc{T}$, and prove that for all proper thick ideals $\mc{I}$ of $\mc{T}$, 
there exist finitely many minimal primes over $\mc{I}$ and the product of their 1-morphism sets (with repetitions) is contained in $\mc{I}_1$. 
\subsection{Noetherian abelian 2-categories}
\label{noeth}
\bde{noetherian} (1) An abelian 2-category will be called {\em{left (resp. right) noetherian}} if it satisfies the ascending chain condition on thick left (resp. right) ideals. 

(2) An abelian 2-category will be called {\em{noetherian}} if it is both left and right noetherian.

(3) An abelian 2-category will be called {\em{weakly noetherian}} if it satisfies the ascending chain condition on (two-sided) thick ideals. 
\ede
More concretely, an abelian 2-category is noetherian if for every chain of thick left ideals 
\[
\mc{I} \subseteq \mc{I}_2 \subseteq \ldots
\] 
there exists an integer $k$ such that $\mc{I}_k = \mc{I}_{k+1} = \ldots$ and such a property is also satisfied for ascending chains of thick right ideals.
\subsection{Existence of minimal primes}
\label{exists-min-prim}
\ble{min-interm} In any abelian 2-category $\mc{T}$, for every thick ideal $\mc{I}$ and every prime ideal $\mc{P}$ containing $\mc{I}$, there is a minimal prime $\mc{P}'$ such that 
\[
\mc{I} \subseteq \mc{P}' \subseteq \mc{P}.
\]
\ele

\begin{proof}
Let $\chi$ denote the set of primes which contain $\mc{I}$ and are contained in $\mc{P}$. We will use Zorn's lemma to produce a minimal element of this set. 
We first show that any nonempty chain in $\chi$ has a lower bound in $\chi$. Take a nonempty chain of prime ideals in $\chi$, say 
\[
\mc{P}^{(1)} \supseteq \mc{P}^{(2)} \supseteq \ldots
\]
Then define $\mc{Q}= \cap_{i=1}^\infty \mc{P}^{(i)}.$ Since each $\mc{P}^{(i)}$ contains $\mc{I}$ and is contained in $\mc{P}$, $\mc{Q}$ is a thick ideal which has these properties. 
It remains to show that $\mc{Q}$ is a prime ideal. Take $f, g \in \mc{T}_1$ such that $f \circ \mc{T}_1 \circ g \subset \mc{Q}_1$, and $f \not \in \mc{Q}_1$. Then $f$ is not in some 
$\mc{P}^{(i)}_1$. Therefore, $f \notin \mc{P}^{(j)}_1$ for $j \geq i$, and by the primeness of $\mc{P}^{(j)}$,  $g \in \mc{P}^{(j)}_1$, for $j \geq i$ as well.  Therefore, 
$g \in \mc{P}^{(k)}_1$ for all $k$, and thus, $g \in \mc{Q}_1$. This implies that $\mc{Q}$ is prime, and Zorn's lemma completes the proof. 
\end{proof}
\subsection{Finiteness and product properties of minimal primes}
\label{pro-min-prim}
\bth{min-primes} In a weakly noetherian abelian 2-category $\mc{T}$, for every proper thick ideal $\mc{I}$, there exist finitely many minimal prime ideals over $\mc{I}$.
Furthermore, there exists a finite list of minimal prime ideals over $\mc{I}$ (potentially with repetition) $\mc{P}^{(1)},..., \mc{P}^{(m)}$ such that the product 
\[
\mc{P}^{(1)}_1 \circ ...\circ \mc{P}^{(m)}_1 \subseteq \mc{I}_1.
\]
\eth
\begin{proof} Denote the set 
\begin{align*}
\chi := \{ \mc{I} \; \; \mbox{a proper thick ideal of} \; \; \mc{T} \mid 
&\; \nexists \; \mbox{prime ideals} \; \mc{P}^{(1)}, \ldots, \mc{P}^{(m)} \supseteq \mc{I} \\
&\; \mbox{such that} \; \; \mc{P}^{(1)}_1 \circ ...\circ \mc{P}^{(m)}_1 \subseteq \mc{I}_1\}.
\end{align*}
Suppose that $\chi$ is nonempty. By the weakly noetherian property of $\mc{T}$, there exists a maximal element of $\chi$ 
(because every ascending chain in $\chi$ eventually stabilizes). 
Let $\mc{I}$ be a maximal element of $\chi$. The ideal $\mc{I}$ cannot be prime, since $\mc{I} \in \chi$. By \thref{prime2} there exist proper thick ideals $\mc{J}$ and $\mc{K}$ 
such that 
\[
\mc{J}_1 \circ \mc{K}_1 \subset \mc{I}_1,
\] 
where $\mc{J}$ and $\mc{K}$ both properly contain $\mc{I}$. The latter property of $\mc{J}$ and $\mc{K}$, and the maximality of $\mc{I}$, imply 
that $\mc{J}, \mc{K} \notin \chi$. Hence, there exist two collection of prime ideals $\mc{P}^{(1)}, \ldots, \mc{P}^{(m)} \supseteq \mc{J}$ and 
$\mc{Q}^{(1)}, \ldots , \mc{Q}^{(n)} \supseteq \mc{K}$ such that
\[
\mc{P}^{(1)}_1 \circ \ldots \circ \mc{P}^{(m)}_1 \subseteq \mc{J}_1 \quad \mbox{and} \quad
\mc{Q}^{(1)}_1 \circ \ldots \circ \mc{Q}^{(n)}_1 \subseteq \mc{K}_1.
\]
Then 
\[
\mc{P}^{(1)}_1 \circ \ldots \circ \mc{P}^{(m)}_1 \circ \mc{Q}^{(1)}_1 \circ \ldots \circ \mc{Q}^{(n)}_1 \subset \mc{I}_1,
\]
giving a contradiction, since the ideals $\mc{P}^{(i)}$ and $\mc{Q}^{(j)}$ are prime and contain $\mc{I}$. 

Hence, $\chi$ is empty. In other words, for every proper thick ideal $\mc{I}$ of $\mc{T}$ there exist prime ideals
$\mc{P}^{(1)}, \ldots, \mc{P}^{(m)} \supseteq \mc{I}$
such that
\begin{equation}
\label{pr}
\mc{P}^{(1)}_1 \circ \ldots \circ \mc{P}^{(m)}_1 \subseteq \mc{I}_1.
\end{equation}
Applying \leref{min-interm}, we obtain that for each $\mc{P}^{(i)}$, there exists a minimal prime $\overline{\mc{P}}^{(i)}$ 
over $\mc{I}$ such that $\overline{\mc{P}}^{(i)} \subseteq \mc{P}^{(i)}$. Combining this with \eqref{pr} gives
\[
\overline{\mc{P}}^{(1)}_1 \circ \ldots \circ \overline{\mc{P}}^{(m)}_1 \subseteq \mc{I}_1
\]
for the minimal primes $\overline{\mc{P}}^{(1)}, \ldots, \overline{\mc{P}}^{(m)}$ of $\mc{I}$. 

Finally, we claim that every minimal 
prime ideal $\mc{P}$ over $\mc{I}$ is in the list $\overline{\mc{P}}^{(1)}, \ldots, \overline{\mc{P}}^{(m)}$. This implies that there are only finitely many 
primes of $\mc{T}$ that are minimal over $\mc{I}$. Indeed, we have 
\[
\overline{\mc{P}}^{(1)}_1 \circ \ldots \circ \overline{\mc{P}}^{(m)}_1 \subseteq \mc{P}_1,
\]
and by the primeness of $\mc{P}$, we have 
\[
\mc{I} \subseteq \mc{P}^{(i)} \subseteq \mc{P}
\]
for some $i$. Since $\mc{P}$ is minimal over $\mc{I}$, $\mc{P}^{(i)} = \mc{P}$.
\end{proof}

The following corollary follows from applying \thref{min-primes} to $0_{\mc{T}}$. 

\bco{fin-noeth} A weakly noetherian abelian 2-category has finitely many minimal prime ideals. 
\eco 

We also have the following corollary of \thref{min-primes}:
\bco{Zaris-w-nothe} For a weakly noetherian abelian 2-category $\mc{T}$, all closed subsets of $\Spec (\mc{T})$ (with respect to the Zariski topology) 
are finite intersections of finitely many sets of the form $V(\mc{P})$ for prime ideals $\mc{P}$ of $\mc{T}$, recall \S \ref{Zariski}.
\eco
\sectionnew{The completely prime spectrum and the semiprime spectrum}
\label{5}
In this section we define the notions of completely prime and semiprime ideals of abelian 2-categories, and give equivalent characterizations, one 
of which is an extension of the Levitzki--Nagata theorem from noncommutative ring theory. 
\subsection{Completely prime ideals}
\bde{compl}
A thick ideal $\mc{P}$ of an abelian 2-category $\mc{T}$ will be called {\em{completely prime}} when it has the property that 
for all $f, g \in \mc{T}_1$:
\[
f \circ g \subseteq \mc{P}_1 \quad \Rightarrow \quad f \in \mc{P}_1 \; \; \mbox{or} \; \; g \in \mc{P}_1.
\]
\ede
This is equivalent to saying that for all 1-morphisms $f$ and $g$ of $\mc{T}$, if $fg$ is not defined or $fg$ is a 1-morphism in $\mc{P}$, then 
$f$ or $g$ is a 1-morphism in $\mc{P}$.
The stronger assumption, including the case of the condition when $fg$ is not defined, is needed to get the correct 
analog of a completely prime ideal of an algebra with a set of orthogonal idempotents. Let $R$ be a ring and $\{e_s\}$ 
be a collection of orthogonal idempotents. If $I$ is a completely prime ideal of $R$ such that 
\[
I \subseteq \bigoplus_{s,t} e_s R e_t,
\]
then for all $s, t \neq t', s'$,
\[
\mbox{either} \quad  e_s R e_t \subseteq I \quad \mbox{or} \quad e_{t'} R e_{s'} \subseteq I,
\]
because $(e_s R e_t) (e_{t'} R e_{s'}) =0$. 

\thref{prime1} implies the following:
\bco{cpr}
Every completely prime ideal of an abelian 2-category is prime.
\eco
\begin{proof} Assume that $\mc{P}$ is a completely prime ideal of $\mc{T}$. Let $f \in \mc{T}(A_3,A_4)$ and $g \in \mc{T}(A_1,A_2)$ be such that
\[
f \circ \mc{T}_1 \circ g \subseteq \mc{P}_1.
\]
If $A_2 \neq A_3$, then $fg$ is not defined and the assumption on $\mc{P}$ gives that either $f \in \mc{P}_1$ or $g \in \mc{P}_1$. If $A_2 = A_3$, then
\[
fg = f 1_{A_2} g \in f \circ \mc{T}_1 \circ g \subseteq \mc{P}_1,
\]
and, again by the assumption on $\mc{P}$, we have that either $f \in \mc{P}_1$ or $g \in \mc{P}_1$. 
\end{proof}

For every abelian 2-category $\mc{T}$, given an object $A$ of $\mc{T}$, consider the 2-subcategory $\mc{T}_A$ of $\mc{T}$ having one object $A$ 
and such that $\mc{T}_{A}(A,A) := \mc{T}(A,A)$. It is an abelian 2-category with one object (i.e., a multiring category). The next lemma
shows that the completely prime ideals of an abelian 2-category $\mc{T}$ are classified in terms of the completely prime ideals of these 
multiring categories.

\ble{c-prime} Let $\mc{T}$ be an abelian 2-category.
\begin{enumerate}
\item If $\mc{P}$ is a completely prime ideal of $\mc{T}$, then there exists an object $A$ of $\mc{T}$ and 
a completely prime ideal $\mc{Q}$ of the multiring category $\mc{T}_A$ such that
\begin{equation}
\label{PQ}
\mc{P}(B,C) =
\begin{cases}
\mc{Q}(A,A), &\mbox{if} \; \; B=C =A
\\
\mc{T}(B,C), &\mbox{otherwise}.
\end{cases}
\end{equation}
\item If $A$ is an object of $\mc{T}$ and $\mc{Q}$ is a completely prime ideal of $\mc{T}_A$ such that 
\begin{equation}
\label{PQcond}
\mc{T}(B,A) \circ \mc{T}(A,B) \subseteq \mc{Q}(A,A)
\end{equation}
for every object $B$ of $\mc{T}$, then \eqref{PQ} defines a completely prime ideal $\mc{P}$ of $\mc{T}$.
\end{enumerate}
\ele
\begin{proof} (1) Since $\mc{P}$ is a proper thick ideal of $\mc{T}$, there exists an object $A$ of $\mc{T}$ such that $1_A \notin \mc{P}_1$. 
(Otherwise $\mc{P}_1$ will contain all 1-morphisms of $\mc{T}$ because
\begin{equation}
\label{1coT}
\mc{T}(B,A) \circ 1_A = \mc{T}(B,A).
\end{equation}
This will contradict the properness of $\mc{P}$.) Obviously 
\[
\mc{Q}(A,A) := \mc{P}(A,A)
\]
defines a completely prime ideal of the multiring category $\mc{T}_A$. It remains to show that $\mc{P}$ 
is given by \eqref{PQ} in terms of $\mc{Q}$. 

If $B$ is an object of $\mc{T}$ which is different from $A$, then the composition $1_A 1_B$ is not defined 
and $1_A \notin \mc{P}_1$, hence $1_B \in \mc{P}_1$. It follows that $\mc{P}$ is given by \eqref{PQ} by an argument 
similar to \eqref{1coT}.

(2) The condition \eqref{PQcond} ensures that the weak thick subcategory $\mc{P}$ of $\mc{T}$ given by \eqref{PQ}, is a thick ideal of $\mc{T}$. Its complete 
primeness is easy to show. 
\end{proof}
\bde{domain} A multiring category $\mc{T}$ will be called a {\em{domain}} if its zero ideal $0_\mc{T}$ is completely prime, i.e., if
\[
M \otimes N \cong 0 \quad \Rightarrow \quad M \cong 0 \; \; \mbox{or} \; \; N \cong 0
\]
for all objects $M$ of $\mc{T}$. 

An abelian 2-category $\mc{T}$ will be called {\em{prime}}, it its zero ideal $0_{\mc{T}}$ is prime.
\ede
\bex{H-mod} Let $H$ be a Hopf algebra over a field $\kk$. Denote by $H\mathrm{-mod}$ the category of finite dimensional $H$-modules. 
It is a $\kk$-linear multiring category. This category is a domain: if $V, W \in H\mathrm{-mod}$ are such that $V \otimes W \cong 0$, then
\[
\dim V \dim W =0.
\]
Therefore, either $\dim V=0$ or $\dim W=0$. So, either $V \cong 0$ or $W \cong 0$. 
\qed
\eex
Let $\mc{T}$ be an abelian 2-category. The same proof shows that if
\begin{enumerate}
\item $\eta : \mc{T}_1 \to R$ is a map such that $R$ is a domain and $\eta(fg) = \eta(f) \eta(g)$ for all $f, g \in \mc{T}$ for which the composition is defined, and
\item $\mc{I}$ is a thick ideal of $\mc{T}$ such that $\mc{I}_1 = \eta^{-1}(0)$, 
\end{enumerate}
then $\mc{I}$ is a completely prime ideal of $\mc{T}$. 
\subsection{Semiprime ideals}
\label{semip}
\bde{semip} A thick ideal of an abelian 2-category will be called {\em{semiprime}} if it is an intersection of prime ideals. 
An abelian 2-category $\mc{T}$ will be called {\em{semiprime}}, it its zero ideal $0_{\mc{T}}$ is semiprime.
\ede
\thref{min-primes} implies that in a weak noetherian abelian 2-category every semiprime ideal is the intersection of the finitely many 
minimal primes over it.

The following theorem is a categorical version of the Levitzki--Nagata theorem.

\bth{semi-p}
A thick ideal $\mc{Q}$ is semiprime if and only if for all $f \in \mc{T}_1$, 
\begin{equation}
\label{cond}
f \circ \mc{T}_1 \circ f \subseteq \mc{Q}_1 \quad \Rightarrow \quad f \in \mc{Q}_1.
\end{equation}
\eth
\begin{proof} First, suppose $\mc{Q}=\bigcap_{s} \mc{P}^{(s)}$ for some collection $\{\mc{P}^{(s)}\}$ of primes of $\mc{T}$. 
Suppose $f \in \mc{T}_1$, and $f \circ \mc{T}_1 \circ f \subseteq \mc{Q}_1.$ 
By primeness, $f \in \mc{P}^{(s)}_1$ for all $s$. Therefore, $f \in \mc{Q}_1$.

For the other direction, suppose that $\mc{Q}$ is a thick ideal of $\mc{T}$ having the property \eqref{cond}.
Choose an element 
\[
g \in \mc{T}_1, \; \; g \notin \mc{Q}_1,
\]
and set $g_0:=g$. It follows from \eqref{cond} that $g_0 \circ \mc{T}_1 \circ g_0 \not \subseteq \mc{Q}_1$. 
Choose
\[
g_1 \in g_0 \circ \mc{T}_1 \circ g_0, \; \; g_1 \notin \mc{Q}_1. 
\]
Again, since $g_1 \notin \mc{Q}_1$, the condition \eqref{cond} implies that $g_1 \circ \mc{T}_1 \circ g_1 \not \subseteq \mc{Q}_1$. Proceeding inductively in this manner, 
we construct a sequence of 1-morphisms $g_0, g_1, \ldots$ of $\mc{T}$ such that  
\begin{equation}
\label{gi}
g_i \in g_{i-1} \circ \mc{T}_1 \circ g_{i-1}, \; \; g_i \notin \mc{Q}_1.
\end{equation}
Since $g_i \in g_{i-1} \circ \mc{T} \circ g_{i-1}$, we have $g_i \circ \mc{T} \circ g_i \subseteq g_{i-1} \circ \mc{T} \circ g_{i-1}.$ Consider the set $S$ of thick ideals $\mc{I}$ 
of $\mc{T}$ such that 
\[
\mc{Q} \subseteq \mc{I} \quad \mbox{and} \quad
g_i \notin \mc{I}_1 \quad \mbox{for all} \quad i=0,1, \ldots.
\]
This set is nonempty because $\mc{Q} \in S$. Since the union of a chain of thick ideals is a thick ideal, we can apply Zorn's lemma to 
get that $S$ contains a maximal element. Denote one such element by $\mc{P}$. The proper thick ideal $\mc{P}$ is prime. Indeed, if $\mc{J}$ and $\mc{K}$ 
are thick ideals that properly contain $\mc{P}$, then by maximality of $\mc{P}$, there are some $g_j \in \mc{J}_1$ and $g_k \in \mc{K}_1$. 
If $m$ is the max of $j$ and $k$, then $g_m$ is in both $\mc{J}_1$ and $\mc{K}_1$ by the first property in \eqref{gi}. Hence,
\[
g_{m+1} \in g_m \circ \mc{T}_1 \circ g_m \subseteq \mc{J}_1 \circ \mc{K}_1 \quad \mbox{and} \quad 
g_{m+1} \not \in \mc{P}_1.
\]
Therefore, $\mc{J}_1 \circ \mc{K}_1 \not \subseteq \mc{P}_1$, and by \thref{prime2}, $\mc{P}$ is prime. For every element $g \in \mc{T}_1$ that is not in 
$\mc{Q}_1$, we have produced a prime $\mc{P}^{(g)}$ of $\mc{T}$ such that
\[
\mc{Q} \subseteq \mc{P}^{(g)} \quad \mbox{and} \quad g \notin \mc{P}^{(g)}. 
\]
Therefore, 
\[
\mc{Q} = \bigcap_{g \in \mc{T}_1 \backslash \mc{Q}_1} \mc{P}^{(g)},
\]
which completes the proof of the theorem.
\end{proof}

\bth{semi-pr-char} Suppose $\mc{Q}$ is a proper thick ideal in an abelian 2-category $\mc{T}$. Then the following are equivalent:
\begin{enumerate}
\item $\mc{Q}$ is semiprime;
\item If $\mc{I}$ is any thick ideal of $\mc{T}$ such that $\mc{I}_1 \circ \mc{I}_1 \subseteq \mc{Q}_1$, then $\mc{I} \subseteq \mc{Q}$;
\item If $\mc{I}$ is any thick ideal properly containing $\mc{Q}$, then $\mc{I}_1 \circ \mc{I}_1\not \subseteq \mc{Q}_1$;
\item If $\mc{I}$ is any right thick ideal of $\mc{T}$ such that $\mc{I}_1 \circ \mc{I}_1 \subseteq \mc{Q}_1$, then $\mc{I} \subseteq \mc{Q}$;
\item If $\mc{I}$ is any left thick ideal of $\mc{T}$ such that $\mc{I}_1 \circ \mc{I}_1 \subseteq \mc{Q}_1$, then $\mc{I} \subseteq \mc{Q}$.
\end{enumerate}
\eth
\begin{proof}
(1) $\Rightarrow$ (4): Suppose $\mc{Q}$ is semiprime and $\mc{I}$ is a right thick ideal with $\mc{I}_1 \circ \mc{I}_1 \subset \mc{Q}_1$. 
Take any $i \in \mc{I}_1$. Then $i \circ t \in \mc{I}$ for all $t \in \mc{T}_1$. Therefore, $i \circ \mc{T}_1 \circ i \in \mc{Q}_1$. 
\thref{semi-p} implies that $i \in \mc{Q}$. Hence, $\mc{I}_1 \subseteq \mc{Q}_1$, and thus $\mc{I} \subseteq \mc{Q}$ by \reref{incl-thick}.

(1) $\Rightarrow$ (5): This follows from a symmetric argument.

(4) $\Rightarrow$ (5) and (3): This is clear, since a thick ideal is also a right thick ideal, and a left thick ideal.

(3) $\Rightarrow$ (2): Suppose (3) holds, and $\mc{I}$ is a thick ideal with $\mc{I}_1 \circ \mc{I}_1 \subseteq \mc{Q}_1$. 
Then $\langle \mc{I}_1 \cup \mc{Q}_1 \rangle$ is a thick ideal containing $\mc{Q}$. 
Since 
\[
(\mc{I}_1 \cup \mc{Q}_1) \circ (\mc{I}_1 \cup \mc{Q}_1)= (\mc{I}_1 \circ \mc{I}_1) \cup (\mc{Q}_1 \circ \mc{I}_1) \cup (\mc{I}_1 \circ \mc{Q}_1) \cup (\mc{Q}_1 \circ \mc{Q}_1) 
\subseteq \mc{Q}_1, 
\]
applying \leref{MTN}, we obtain
\[
\langle \mc{I}_1 \cup \mc{Q}_1 \rangle_1 \circ \langle \mc{I}_1 \cup \mc{Q}_1 \rangle_1 \subseteq 
\langle (\mc{I}_1 \cup \mc{Q}_1) \circ \mc{T}_1 \circ (\mc{I}_1 \cup \mc{Q}_1) \rangle_1
= \langle (\mc{I}_1 \cup \mc{Q}_1) \circ (\mc{I}_1 \cup \mc{Q}_1) \rangle_1 \subseteq \mc{Q}_1.
\]
From the assumption that the ideal $\mc{Q}$ has the property (3) and the fact that $\langle \mc{I}_1 \cup \mc{Q}_1 \rangle$ is a thick ideal containing $\mc{Q}$, 
we get that $\langle \mc{I}_1 \cup\mc{Q}_1 \rangle = \mc{Q}$. Therefore, $\mc{I}_1 \subseteq \mc{Q}_1$, and thus  $\mc{I} \subseteq \mc{Q}$ by \reref{incl-thick}. 

(2) $\Rightarrow$ (1): Suppose (2) holds, and $f \in \mc{T}$ is a 1-morphism such that $f \circ \mc{T}_1 \circ f \subseteq \mc{Q}_1$. \leref{MTN} implies that
\[
\langle f \rangle_1 \circ \langle f \rangle_1 \subseteq \langle f \circ \mc{T}_1 \circ f \rangle_1 \subseteq \mc{Q}_1.
\]
Therefore, by (2), $\langle f \rangle_1 \subseteq \mc{Q}_1$, and so, $f \in \mc{Q}_1$. Hence, $\mc{Q}_1$ is semiprime.
\end{proof}
We have the following corollary from the characterizations (4) and (5) of semiprime ideals in the previous theorem. For a subset $X \subseteq \mc{T}_1$, 
denote by $X^{\circ n} := X \circ \cdots \circ X$ the $n$-fold composition power.
\ble{powers-semipr}
If $\mc{Q}$ is a semiprime ideal of the abelian 2-category $\mc{T}$, 
and $\mc{I}$ is a right or left thick ideal with $(\mc{I}_1)^{\circ n} \subseteq \mc{Q}_1$, then $\mc{I} \subseteq \mc{Q}$. 
\ele
\begin{proof}
We prove the statement by induction on $n$. For $n \geq 2$, we have
\[
((\mc{I}_1)^{\circ (n-1)})^{\circ 2} = (\mc{I}_1)^{\circ n} \circ (\mc{I}_1)^{n-2} \subseteq \mc{Q}_1.
\]
\thref{semi-pr-char} implies that $(\mc{I}_1)^{\circ (n-1)} \subseteq  \mc{Q}_1$, and so by the inductive assumption, $\mc{I} \subseteq \mc{Q}$.
\end{proof}
\sectionnew{The Serre prime spectra of abelian 2-categories and $\Zset_+$-rings}
\label{6}
In this section we define and investigate the notions of Serre prime, semiprime, and completely prime  ideals of abelian 2-categories 
and $\Zset_+$-rings. We establish that the corresponding topological spaces for abelian 2-categories and $\Zset_+$-rings are homeomorphic.
We also describe the relations of the first set of notions to the notions of prime, completely prime and semiprime ideals of abelian 2-categories, 
and the second set of notions to the prime spectra of rings.
\subsection{Serre ideals of abelian 2-categories}
\label{Serre}
Recall that a Serre subcategory of an abelian 1-category is a subcategory which is closed under subobjects, quotients, and extensions. 
Every Serre subcategory $\mc{I}$ of an abelian category $\mc{C}$ is thick, and in particular, is closed 
under isomorphisms.
For such a subcategory, one forms the Serre quotient $\mc{C}/\mc{I}$ which has a canonical structure of 
abelian category \cite[\S 10.3]{Weibel1}. 
By \cite[Theorem 5]{Quillen1}, for every Serre subcategory $\mc{I}$ of an abelian category $\mc{C}$,
we have the exact sequence
\begin{equation}
\label{Quillen}
K_0(\mc{I}) \to K_0(\mc{C}) \to K_0(\mc{C}/\mc{I}) \to 0.
\end{equation}

\bde{Serre} (1) We call a thick ideal $\mc{I}$ of an abelian 2-category $\mc{T}$ a {\em{Serre ideal}} if for every two objects $A_1, A_2 \in \mc{T}$, 
\[
\mbox{$\mc{I}(A_1,A_2)$ is a Serre subcategory of $\mc{T}(A_1, A_2)$}.
\]

(2) A {\em{Serre prime (resp. semiprime, completely prime) ideal}} $\mc{P}$ of an abelian 2-category $\mc{T}$
is a prime (resp. semiprime, completely prime) ideal which is a Serre ideal.
\ede

In the terminology of \deref{weak-sub}, a Serre ideal of an abelian 2-category $\mc{T}$ is a weak subcategory $\mc{I}$ with the same set of objects such that
\begin{enumerate}
\item for any pair of objects $A_1, A_2 \in \mc{T}$, $\mc{I}(A_1, A_2)$ is a Serre subcategory of the abelian category $\mc{T}(A_1,A_2)$ and
\item $\mc{I}_1 \circ \mc{T}_1 \subseteq \mc{I}_1$, $\mc{T}_1 \circ \mc{I}_1 \subseteq \mc{I}_1$.
\end{enumerate}

We will say that $\mc{I}$ is a {\em{left (resp. right) Serre ideal}} of $\mc{T}$ if condition (1) is satisfied and $\mc{T}_1 \circ \mc{I}_1 \subseteq \mc{I}_1$
(resp. $\mc{I}_1 \circ \mc{T}_1 \subseteq \mc{I}_1$).

\bpr{Serre} For every Serre ideal $\mc{I}$ of an abelian 2-category $\mc{T}$ such that $1_A \notin \mc{I}(A,A)$ for all objects $A \in \mc{T}$, one can form 
the Serre quotient $\mc{T}/\mc{I}$ with the same set of objects, with the morphism 1-categories 
\[
(\mc{T}/\mc{I})(A_1, A_2) := \mc{T}(A_1, A_2)/\mc{I}(A_1, A_2) \quad \mbox{for} \quad A_1, A_2 \in \mc{T},
\] 
and with identity 1-morphisms given by the images of $1_A$. This quotient is an abelian 2-category.
\epr
The proof of the proposition is direct, using \eqref{Quillen} and the following well known fact:

{\em{If, for $i=1,2$, $\mc{C}_i$ are abelian categories, $\mc{I}_i$ are Serre subcategories, and $F : \mc{C}_1 \to \mc{C}_2$ 
is an exact functor such that $F(\mc{I}_1) \subseteq \mc{I}_2$, then the induced functor $\overline{F} : \mc{C}_1 / \mc{I}_1 \to \mc{C}_2 / \mc{I}_2$ 
is exact.
}}

This follows from the commutativity of the square diagram consisting of the compositions of functors 
$\mc{C}_1 \to \mc{C}_1 / \mc{I}_1 \stackrel{\overline{F}}{\to} \mc{C}_2 / \mc{I}_2$ and $\mc{C}_1\stackrel{F}{\to} \mc{C}_2 \to \mc{C}_2 / \mc{I}_2$, 
the exactness of the projection functors $\mc{C}_i \to \mc{C}_i / \mc{I}_i$ (see \cite[Exercise 10.3.2(4)]{Weibel1}), and 
the fact that every exact sequence in $\mc{C}_1 / \mc{I}_1$ is isomorphic to one coming from an exact sequence in $\mc{C}_1$,
\cite{Gabriel1}.

It is easy to prove that, similarly to the ring theoretic case, we have the following:
\ble{compl-pr} A proper Serre ideal $\mc{P}$ of a multiring category $\mc{T}$ is completely prime, if and only if 
the Serre quotient $\mc{T}/\mc{I}$ is a domain in the sense of \deref{domain}.
\ele

Analogously to Lemmas \ref{lint-thick} and \ref{lMTN} one proves the following result. We leave the details to the reader.
\ble{Serre} Let $\mc{T}$ be an abelian 2-category.

(1) The intersection of any family of Serre ideals of $\mc{T}$ is a Serre ideal of $\mc{T}$. 
In particular, for any subset $\mc{M} \subseteq \mc{T}_1$, there exists a unique minimal Serre ideal of $\mc{T}$ containing $\mc{M}$; 
it will be denoted by $\langle \mc{M} \rangle^S$. 

(2) For $\mc{M}, \mc{N} \subseteq \mc{T}_1$, we have
\[
\langle \mc{M} \rangle^S_1 \circ \langle \mc{N} \rangle^S_1 \subseteq \langle \mc{M} \circ \mc{T}_1 \circ \mc{N} \rangle^S_1.
\]
\ele
\subsection{Serre prime ideals of abelian 2-categories}
\label{Serre-prime}
Similarly to the proofs of Theorems \ref{tprime1}, \ref{tprime2}, \ref{tmult-max} and \ref{tsemi-pr-char}, using \leref{Serre}, one proves the following
result:
\bth{Serre-pr} Let $\mc{T}$ be an abelian 2-category.
\begin{enumerate}
\item The following are equivalent for a proper Serre ideal $\mc{P}$ of $\mc{T}$:
\begin{itemize}
\item[(a)] $\mc{P}$ is a Serre prime ideal;
\item[(b)] If $\mc{I}$ and $\mc{J}$ are any Serre ideals of $\mc{T}$ such that $\mc{I}_1 \circ \mc{J}_1 \subseteq \mc{P}_1$, then either $\mc{I} \subseteq \mc{P}$
or  $\mc{J} \subseteq \mc{P}$;
\item[(c)] If $\mc{I}$ and $\mc{J}$ are any Serre ideals properly containing $\mc{P}$, then $\mc{I}_1 \circ \mc{J}_1\not \subseteq \mc{P}_1$;
\item[(d)] If $\mc{I}$ and $\mc{J}$ are any left Serre ideals of $\mc{T}$ such that $\mc{I}_1 \circ \mc{J}_1 \subseteq \mc{P}_1$, then either 
$\mc{I} \subseteq \mc{P}$ or $\mc{J} \subseteq \mc{P}$. 
\end{itemize}
\item Let $\mc{M}$ be a nonempty multiplicative subset of $\mc{T}_1$ (cf. \deref{mult}) and 
$\mc{I}$ be a Serre ideal of $\mc{T}$ such that $\mc{I}_1 \cap \mc{M} = \varnothing$.
Let $\mc{P}$ be a maximal element of the collection of Serre ideals of $\mc{T}$ containing $\mc{I}$ and intersecting 
$\mc{M}$ trivially, equipped with the inclusion relation, i.e., $\mc{P}$ is a maximal element of the set 
\[
X(\mc{M}, \mc{I}) := \{ \mc{K} \; \mbox{a Serre ideal of} \; \; \mc{T} \mid \mc{K} \supseteq \mc{I} \; \; \mbox{and} \; \; \mc{K}_1 \cap \mc{M} = \varnothing \}.
\]
Then $\mc{P}$ is Serre prime ideal.
\item The following are equivalent for a proper Serre ideal $\mc{Q}$ of $\mc{T}$:
\begin{itemize}
\item[(a)] $\mc{Q}$ is a Serre semiprime ideal;
\item[(b)] If $\mc{I}$ is any Serre ideal of $R$ such that $\mc{I}_1 \circ \mc{I}_1 \subseteq \mc{Q}_1$, then $\mc{I} \subseteq \mc{Q}$;
\item[(c)] If $\mc{I}$ is any Serre ideal properly containing $\mc{Q}$, then $\mc{I}_1 \circ \mc{I}_1\not \subseteq \mc{Q}_1$;
\item[(d)] If $\mc{I}$ is any left Serre ideal of $\mc{T}$ such that $\mc{I}_1 \circ \mc{I}_1 \subseteq \mc{Q}_1$, then $\mc{I} \subseteq \mc{Q}$.
\end{itemize}
\end{enumerate}
\eth
In the proof of part (1) of the theorem, the key step is to show that a proper Serre ideal $\mc{I}$ of $\mc{T}$ satisfying the property (b) is a Serre prime ideal of $\mc{T}$. 
This is proved by showing that property (b) implies that for all $m, n \in \mc{T}_1$, 
\[
m \circ \mc{T}_1 \circ n \subseteq \mc{P}_1 \quad \Rightarrow m \in \mc{P}_1 \; \; \mbox{or}  \; \; n \in \mc{P}_1.
\]
This fact is verified by repeating the proof of \thref{prime1}, but using \leref{Serre}(2) in place of \leref{MTN}.

The set $X(\mc{M}, \mc{I})$ in part (3) of the theorem is nonempty because $\mc{I} \in X(\mc{M}, \mc{I})$. 
The union of an ascending chain of Serre ideals in the set $X(\mc{M}, \mc{I})$ is obviously a Serre ideal of $\mc{T}$. 
By Zorn's lemma, the set $X(\mc{M}, \mc{I})$ always contains at least one maximal element. 

Similarly to the proof of \coref{nonempty}, we obtain:
\bco{nonempty2} For every proper Serre ideal $\mc{I}$ of an abelian 2-category $\mc{T}$, there 
exists a Serre prime ideal $\mc{P}$ of $\mc{T}$ that contains $\mc{I}$. 
\eco

Analogously to the proof of \thref{min-primes} one proves the following:

\bpr{min-primes} For every abelian 2-category $\mc{T}$ satisfying the ACC on (2-sided) Serre ideals, given a proper Serre ideal $\mc{I}$ of $\mc{T}$, 
there exist finitely many minimal Serre prime ideals over $\mc{I}$. 
Furthermore, there is a finite list of minimal Serre prime ideals over $\mc{I}$ (possibly with repetition) $\mc{P}^{(1)},..., \mc{P}^{(m)}$ such that the product 
\[
\mc{P}^{(1)}_1 \circ ...\circ \mc{P}^{(m)}_1 \subseteq \mc{I}_1.
\]
\epr

Let $\mathrm{Serre}\mbox{-}\mathrm{Spec} (\mc{T})$ denote the set of Serre prime ideals of an abelian 2-category $\mc{T}$. Similarly to 
\S \ref{Zariski}, one shows that it is a topological space with closed subsets given by 
\[
V^S(\mc{I}) = \{ \mc{P} \in \mathrm{Serre}\mbox{-}\mathrm{Spec} (\mc{T}) \mid \mc{P}  \supseteq \mc{I} \}
\] 
for the Serre ideals $\mc{I}$ of $\mc{T}$.
We will refer to this as to the Zariski topology of $\mathrm{Serre}\mbox{-}\mathrm{Spec} (\mc{T})$.
\prref{min-primes} implies that, if $\mc{T}$ satisfies the ACC on Serre ideals, 
then every closed subset of $\mathrm{Serre}\mbox{-}\mathrm{Spec} (\mc{T})$ is a finite intersection of subsets of the 
form $V^S(\mc{P})$ for some $\mc{P} \in \mathrm{Serre}\mbox{-}\mathrm{Spec} (\mc{T})$. In particular, this property holds 
for all weakly noetherian abelian 2-categories $\mc{T}$.

The set-theoretic inclusion 
\begin{equation}
\label{embed-S}
\mathrm{Serre}\mbox{-}\mathrm{Spec} (\mc{T}) \hookrightarrow \mathrm{Spec} (\mc{T})
\end{equation}
realizes $\mathrm{Serre}\mbox{-}\mathrm{Spec} (\mc{T})$ as a topological subspace of $\mathrm{Spec} (\mc{T})$. Indeed, \leref{Serre}(1) implies that 
for every thick ideal $\mc{I}$ of $\mc{T}$ we have
\[
V(\mc{I}) \cap \mathrm{Serre}\mbox{-}\mathrm{Spec} (\mc{T}) = V^S(\langle \mc{I} \rangle^S).
\]
\subsection{The Serre prime spectrum as a ringed space}
\label{ringed-space} For the following subsection, assume that $\mc{C}$ is an abelian monoidal category. In the case when it is strict this is the same as an abelian 2-category with one object. All constructions and results in the paper are valid for abelian monoidal categories without the strictness assumption.
By \reref{different} and the embedding \eqref{embed-S}, 
the Zariski topology we have thus far endowed $\mathrm{Serre}\mbox{-}\mathrm{Spec}(\mc{C})$ with is different from the topology used by Balmer in \cite{Balmer1}. 
The motivation for this consists of the applications to categorification, which we develop below. However, if $\mc{C}$ is an abelian monoidal category, we can 
consider an analogue of Balmer's topology on $\mathrm{Serre}\mbox{-}\mathrm{Spec}(\mc{C}),$ where we define the closed sets of $\mathrm{Serre}\mbox{-}\mathrm{Spec}(\mc{C})$ to be 
\[
V_B^S(X)=\{ \mc{P} \in \mathrm{Serre}\mbox{-}\mathrm{Spec}(\mc{C}) \mid X \cap \mc{P} = \emptyset\}
\]
for any set of objects $X$ in $\mc{C}$. Analogously to Section 2 of \cite{Balmer1}, one shows that this collection defines a topological space. It may be equipped with a sheaf of commutative rings in a similar manner to \cite{Balmer1}. 
Let $U=V^c$ be an open set of $\mathrm{Serre}\mbox{-}\mathrm{Spec}(\mc{C})$, 
where $V=V_B^S(X)$ for some family of objects $X$. 
Let $\mc{C}_V=\bigcap_{\mc{P} \in U} \mc{P}.$ Note that $\mc{C}_V$ is a Serre ideal, since it is an intersection of Serre ideals. We define a presheaf of commutative rings in the following way:
\[
U \mapsto \End_{\mc{C}/ \mc{C}_V}(\overline{1}, \overline{1}),
\]
where $\overline{1}$ is the image of 1 (the unit object of $\mc{C}$ with respect to the monoidal product) in the Serre quotient $\mc{C}/\mc{C}_V$. Recalling \prref{Serre}, $\mc{C}/\mc{C}_V$ has a canonical structure as an abelian monoidal category. By e.g. Proposition 2.2.10 in \cite{Etingof1}, $\End_{\mc{C}/\mc{C}_V}(\overline{1}, \overline{1})$ is a commutative ring. The sheafification of this presheaf gives $\mathrm{Serre}\mbox{-}\mathrm{Spec}(\mc{C})$ the structure of a ringed space. The question about the construction of a ringed space structure on
the spectra of abelian monoidal categories was raised by Michael Wemyss.
\subsection{$\Zset_+$-rings and their Serre prime ideals}
\label{Z+rings} Recall that $\Zset_+ := \{ 0, 1, \ldots \}$.

We will use the following slightly weaker terminology for $\Zset_+$-rings compared to \cite[Definition 3.1.1]{Etingof1}: 
\bde{Z+rings} We will call a ring $R$ a $\Zset_+$-ring if it is a free abelian group and has a $\Zset$-basis $\{b_\gamma \mid \gamma \in \Gamma\}$ 
such that for all $\alpha, \beta \in \Gamma$,
\[
b_\alpha b_\beta = \sum_{\gamma \in \Gamma} n_{\alpha, \beta}^\gamma b_\gamma
\]
for some $n_{\alpha, \beta}^\gamma \in \Zset_+$. 
\ede
In addition, \cite[Definition 3.1.1]{Etingof1} requires that a $\Zset_+$-ring $R$ be a unital ring and 
\begin{equation}
\label{unit}
1 = \sum_{\gamma \in \Gamma} n_\gamma b_\gamma \quad \mbox{for some} \quad n_\gamma \in \Zset_+.
\end{equation}
We do not require a $\Zset_+$-ring to be unital and to have the above additional property in 
order to apply the notion to the Grothendieck rings of abelian 2-categories with infinitely many objects.

For a $\Zset_+$-ring $R$, denote
\[
R_+ := \bigoplus_{\gamma \in \Gamma} \Zset_+ b_\gamma.
\]
For $r, s \in R$, denote
\[
r \leq s \quad \mbox{if} \quad s-r \in R_+.
\]
\bde{Serre-Z+} Let $R$ be a $\Zset_+$-ring.

(1) A left (resp. right) ideal $I$ of $R$ will be called a {\em{a left (resp. right) Serre ideal}} if it has the properties that 
\begin{equation}
\label{Ser-Z+}
I = (I \cap R_+) - (I \cap R_+) \quad \mbox{and} \quad 
s \in R_+, r \in I \cap R_+, s \leq r \Rightarrow s \in I. 
\end{equation}

(2) A Serre ideal of $R$ is a 2-sided ideal $I$ of $R$ which satisfies \eqref{Ser-Z+}.

(3) A Serre prime ideal of $R$ is a proper Serre ideal $P$ of $R$ that has the property that
\[
I J \subseteq P \quad \Rightarrow \quad I \subseteq P \; \; \mbox{or} \; \; J \subseteq P
\]
for all Serre ideals $I,J$ of $R$.

(4) A Serre semiprime ideal of $R$ is an ideal which is the intersection of Serre prime ideals.

(5) A Serre completely prime ideal of $R$ is a proper Serre ideal $P$ that has the property that 
for all $r, s \in R_+$, 
\[
r s \in P \quad \Rightarrow \quad r \in P \quad \mbox{or} \quad s \in P. 
\]
\ede
For a subgroup $I$ (under addition) of a $\Zset_+$-ring $R$, the property \eqref{Ser-Z+} is equivalent to
\begin{equation}
\label{Serre-form}
I = \bigoplus_{\gamma \in \Gamma'}  \Zset b_\gamma \quad \mbox{for some subset} \quad \Gamma' \subseteq \Gamma.
\end{equation}
In particular, the right and 2-sided Serre ideals of $R$ satisfy \eqref{Serre-form}.
Using this fact one easily proves the following theorem, by following 
the strategy of the proofs of Proposition 3.1, Theorem 3.7 and Corollary 3.8 in \cite{Goodearl1}.
\bth{Serre-Z+} Let $R$ be a $\Zset_+$-ring.
\begin{enumerate}
\item The following are equivalent for a proper Serre ideal $P$ of $R$:
\begin{itemize}
\item[(a)] $P$ is a Serre prime ideal;
\item[(b)] If $I$ and $J$ are two Serre ideals of $R$ properly containing $P$, then $I J \not \subseteq P$;
\item[(c)] If $I$ and $J$ are two left Serre ideals of $R$ such that $IJ \subseteq P$, then either 
$I \subseteq P$ or $J \subseteq P$; 
\item[(d)] For all $\alpha, \beta \in \Gamma$, 
\[
b_\alpha R b_\beta \subseteq P \quad \Rightarrow \quad b_\alpha \in P \; \; \mbox{or} \; \; b_\beta \in P.
\]
\end{itemize}
\item A proper Serre ideal $P$ of $R$ is a completely prime Serre ideal if and only if for all $\alpha, \beta \in \Gamma$,
\[
b_\alpha b_\beta \in P \quad \Rightarrow \quad b_\alpha \in P \; \; \mbox{or} \; \; b_\beta \in P.
\]

\item The following are equivalent for a proper Serre ideal $Q$ of $R$:
\begin{itemize}
\item[(a)] $Q$ is a Serre semiprime ideal;

\item[(b)] If $I$ is any Serre ideal of $R$ such that $I^2 \subseteq Q$, then $I \subseteq Q$;
\item[(c)] If $I$ is any Serre ideal of $R$ properly containing $Q$, then $I^2 \not \subseteq Q$;
\item[(d)] For all $r \in R_+$,
\[
r R r \subseteq P \quad \Rightarrow \quad r \in P.
\] 
\end{itemize}
\end{enumerate}
\eth
Denote by $\mathrm{Serre}\mbox{-}\mathrm{Spec} (R)$ the set of Serre prime ideals of a $\Zset_+$-ring $R$. Similarly to 
\S \ref{Serre-prime}, one proves that $\mathrm{Serre}\mbox{-}\mathrm{Spec} (R)$ is a topological space with closed 
subsets 
\[
V^S(I) = \{ P \in \mathrm{Serre}\mbox{-}\mathrm{Spec} (\mc{T}) \mid P  \supseteq I \}
\] 
for the Serre ideals $I$ of $R$.
We will call this the Zariski topology of $\mathrm{Serre}\mbox{-}\mathrm{Spec} (R)$. 
\subsection{$\Zset_+$-rings and abelian 2-categories}
\label{Zplus} 
For an abelian category $\mc{C}$ denote by $\mc{C}_s$ the equivalence classes of its simple objects.
\ble{fin-ab} Assume that $\mc{C}$ is an abelian category in which every object has finite length. 
Then the following hold:

(1) Every two Jordan-H\"older series of an object of $\mc{C}$ contain the same collections of simple subquotients counted with multiplicities and, as a consequence,
\[
K_0(\mc{C}) \cong \bigoplus_{A \in \mc{C}_s} \Zset [A].
\]

(2) The Serre subcategories of $\mc{C}$ are in bijection with the subsets of $\mc{C}_s$. The Serre subcategory corresponding to a subset $X \subseteq \mc{C}_s$ 
is the full subcategory $\mc{S}(X)$ of $\mc{C}$ whose objects have Jordan-H\"older series with simple subquotients isomorphic to objects in $X$.

(3) For every Serre subcategory $\mc{I}$ of $\mc{C}$, 
\[
K_0(\mc{C}/\mc{I}) \cong K_0(\mc{C})/K_0(\mc{I}).
\]
\ele
\begin{proof}
The first part of the lemma is \cite[Theorem 1.5.4]{Etingof1}. 

(2) Clearly, for every subset $X \subseteq \mc{C}_s$, the subcategory $\mc{S}(X)$ of $\mc{C}$ is Serre. 
Assume that $\mc{I}$ is a Serre subcategory of $\mc{C}$. Denote by $X$ the isomorphism classes of simple objects of $\mc{C}$ which 
belong to $\mc{I}$. Since $\mc{I}$ is closed under taking subquotients and isomorphisms, $\mc{I} \subseteq \mc{S}(X)$. Because
$\mc{I}$ is closed under extensions, $\mc{S}(X) \subseteq \mc{I}$. Thus, $\mc{I} = \mc{S}(X)$.

(3) It easily follows from parts (1) and (2) that the first map in \eqref{Quillen} is injective. The resulting short exact sequence from \eqref{Quillen} 
implies the third part of the lemma. 
\end{proof} 

For an abelian 2-category $\mc{T}$, denote by $(\mc{T}_1)_s$ the isomorphism classes of simple 1-morphisms of $\mc{T}$. 
Recall \deref{weak-sub}.
For a subset $X \subseteq (\mc{T}_1)_s$, denote by $\mc{S}(X)$ the (unique) weak subcategory of $\mc{T}$ such that 
\[
\mc{S}(X)(A_1,A_2):= \mc{S}( X \cap \mc{T}(A_1,A_2))
\]
for all $A_1, A_2 \in \mc{T}$.

\bth{fin-ab} Assume that $\mc{T}$ is an abelian 2-category with the property that every 1-morphism of $\mc{T}$ has finite length.
(In other words, every object of $\mc{T}(A_1, A_2)$ has finite length for all objects $A_1, A_2 \in \mc{T}$.)
Then the following hold:

(1) The weak subcategories $\mc{I}$ of $\mc{T}$ with the property that $\mc{I}(A_1, A_2)$ is a 
Serre subcategory of $\mc{T}(A_1, A_2)$ for all $A_1, A_2 \in \mc{T}$ are parametrized by the subsets of $(\mc{T}_1)_s$.
For $X \subseteq (\mc{T}_1)_s$, the corresponding subcategory is $\mc{S}(X)$. 

(2) The Grothendieck ring $K_0(\mc{T})$ is a $\Zset_+$-ring and
\[
K_0(\mc{T}) \cong \bigoplus_{f \in (\mc{T}_1)_s} \Zset [f].
\] 
If, in addition, $\mc{T}$ has finitely many objects, then $K_0(\mc{T})$ has the property \eqref{unit} and, more precisely, 
\[
1 = \sum_{A \in \mc{T}} [1_A].
\]  

(3) The map $K_0$ defines a bijection between the sets of left (resp. right, 2-sided) Serre ideals of $\mc{T}$ and of $K_0(\mc{T})$.

(4) The map $K_0$ defines a homeomorphism 
\[
K_0 : \mathrm{Serre}\mbox{-}\mathrm{Spec} (\mc{T}) \stackrel{\cong}{\rightarrow} \mathrm{Serre}\mbox{-}\mathrm{Spec} (K_0(\mc{T})).
\]
It is a bijection between the subsets of completely prime (resp. semiprime) ideals of $\mc{T}$ and $K_0(\mc{T})$. 
\eth
\begin{proof} Part (1) follows from \leref{fin-ab}(2).

(2) The fact that $K_0(\mc{T})$ is a $\Zset_+$-ring follows from the fact that for every abelian category $\mc{C}$ and $B \in \mc{C}$,
\[
[B] \in \bigoplus_{A \in \mc{C}_s}  \Zset_+ [A].
\]
The second statement in part (2) is obvious.

(3) We consider the case of left Serre ideals, the other two cases being analogous. Let $\mc{I}$ be a left Serre ideal of $\mc{T}$. 
By part (1) of the theorem, $\mc{I} = \mc{S}(X)$ for some $X \subseteq (\mc{T}_1)_s$. Therefore,
the subset 
\[
K_0(\mc{I}) = \bigoplus_{f \in X} \Zset [f]  \subseteq K_0(\mc{T})
\]
has the property \eqref{Ser-Z+}. Since $\mc{T}_1 \circ \mc{I}_1 \subseteq \mc{I}_1$, we have 
$K_0(\mc{T}) K_0(\mc{I}) \subseteq K_0(\mc{I})$, and thus, $K_0(\mc{I})$ is a left Serre ideal of 
$K_0(\mc{T})$. 

Next, let $I$ be a left Serre ideal of $K_0(\mc{T})$. By \eqref{Serre-form}, 
\[
I = \bigoplus_{f \in X} \Zset [f]
\]
for some $X \subseteq (\mc{T}_1)_s$. Let $\mc{I}$ be the weak subcategory $\mc{S}(X)$ of $\mc{T}$. 
Clearly, $K_0(\mc{I}) =I$. To show that $\mc{I}$ is a left Serre ideal of $\mc{T}$, 
it remains to prove that $\mc{T}_1 \circ \mc{I}_1 \subseteq \mc{I}_1$, i.e., that
\[
g_2 f_1 \in \mc{I}(A_1, A_3) \quad \mbox{for all} 
\quad g_2 \in \mc{T}(A_2, A_3), f_1 \in \mc{I}(A_1, A_2)
\]
for all objects $A_1, A_2, A_3$ of $\mc{T}$. Since $I$ is a left Serre ideal, 
\[
[g_2 f_1] = [g_2] [f_1] \in I =  \bigoplus_{f \in X} \Zset [f],
\]
and thus, $g_2 f_1 \in \mc{S}(X) = \mc{I}$. 

It is straightforward to verify that the above two maps $\mc{I} \mt K_0(\mc{I})$ and $I \mt \mc{I}$ are 
inverse bijections between the left Serre ideals of $\mc{T}$ and $K_0(\mc{T})$.  

(4) Similarly to part (3) one proves that the map $K_0$ defines a bijection between the 
prime (resp. completely prime, semiprime) ideals of the abelian 2-category $\mc{T}$ and the $\Zset_+$-ring $K_0(\mc{T})$. 
In the first case one uses the characterization of Serre prime ideals of an abelian 2-category in  
\thref{Serre-pr}(1)(b) vs. the definition of Serre prime ideals of a $\Zset_+$-ring. 
In the second case one uses the definitions of completely prime ideals in the two settings. 
In the third case one uses the characterizations of Serre semiprime ideals in the two settings 
given in Theorems \ref{tSerre-pr}(3)(b) and \ref{tSerre-Z+}(3)(b).

The fact that the map
\[
K_0 : \mathrm{Serre}\mbox{-}\mathrm{Spec} (\mc{T}) \rightarrow \mathrm{Serre}\mbox{-}\mathrm{Spec} (K_0(\mc{T}))
\]
is a homeomorphism follows from the definitions of the collections of closed sets in the two cases in terms of Serre ideals and the bijection in 
part (3) of the theorem.
\end{proof}

We have the following immediate corollary of part (3) of the theorem and \leref{fin-ab}:
\bco{fact-cat} Let $\mc{T}$ be an abelian 2-category which is a categorification of the $\kk$-algebra $R \otimes_\Zset \kk$ for a 
$\Zset_+$-ring $R$. If $I$ is a Serre ideal of $R$ and $\mc{I}$ is the unique Serre ideal of $\mc{T}$ with $K_0(\mc{I}) = I$ as in 
\thref{fin-ab}(3), then $\mc{T}/\mc{I}$ is a categorification of the $\kk$-algebra $(R/I) \otimes_\Zset \kk$. 
\eco
\subsection{Serre prime ideals of $\Zset_+$-rings vs prime ideals}
\label{Serre-2cat-Z+} 
Let $R$ be a $\Zset_+$-ring. In general, $\mathrm{Serre}\mbox{-}\mathrm{Spec} (R)$ is not a subset of the prime spectrum $\mathrm{Spec}(R)$ of $R$. 
Similarly a Serre completely prime ideal of $R$ is not necessarily a completely prime ideal of $R$ (in the classical sense), and 
a Serre semiprime ideal of $R$ is not necessarily a semiprime ideal of $R$. The point in all three cases is that the notions of Serre type are formulated in terms of 
inclusion properties concerning elements of $R_+$, while the classical notions are formulated in terms of inclusion properties concerning elements of the full ring $R$. 

\bex{first-ex} Consider the commutative $\Zset_+$-ring $R := \Zset[x]/(x^2-1)$ with positive $\Zset$-basis $\{1, x \}$.  The 0-ideal of $R$ is 
Serre prime while it is not a prime ideal of $R$. 
\eex

\bex{second-ex} Consider the setting of \exref{H-mod} and assume that $H$ is a finite dimensional Hopf algebra over the field $\kk$. 
The 0 ideal of $H\mathrm{-mod}$ is Serre completely prime, and by \thref{fin-ab}(4), $0$ is a Serre completely prime ideal of $K_0(H\mathrm{-mod})$.
However, 0 is not a completely prime ideal of the ring $K_0(H\mathrm{-mod})$ because $K_0(H\mathrm{-mod}) \otimes_\Zset \Qset$ 
is a finite dimensional algebra over $\Qset$ and thus definitely has 0 divisors. Furthermore the 0 ideal of $K_0(H\mathrm{-mod})$  
is not even semiprime, except for the special case when the algebra $K_0(H\mathrm{-mod}) \otimes_\Zset \Qset$ is semisimple 
(because the radical of this finite dimensional algebra is nilpotent). 
\qed
\eex

On the other hand, the following lemma provides a simple but important fact about getting Serre prime (resp. completeley prime, semiprime) ideals of a $\Zset_+$-ring $R$ 
from particular types of prime (resp. completeley prime, semiprime) ideals of a $R$ in the classical sense.
\ble{prime-to-Serreprime} Assume that $R$ is a $\Zset_+$-ring with a positive basis $\{ b_\gamma \mid \gamma \in \Gamma \}$. If 
\[
I = \bigoplus_{\gamma \in \Gamma'}  \Zset b_\gamma \quad \mbox{for some subset} \quad \Gamma' \subseteq \Gamma
\]
and $I$ is a prime (resp. completely prime, semiprime) ideal of $R$ in the classical sense, then
$I$ is a Serre prime (resp. completely prime, semiprime) ideals of $R$.
\ele
\begin{proof} The first property of $I$ is equivalent to the one in \eqref{Ser-Z+}. The assumption that $I$ is a prime (resp. completely prime, semiprime) ideal
of $R$ implies that it satisfies the condition (b) in \thref{Serre-Z+}(1) in the first case, the condition in \thref{Serre-Z+}(2) in the second case, and the condition (d) in \thref{Serre-Z+}(3) 
in the third case. For example, if $I$ is a semiprime ideal of $R$ in the classical sense, the it satisfies the condition (d) in \thref{Serre-Z+}(3) for all $r \in R$.
Now the lemma follows from \thref{Serre-Z+}.
\end{proof}
\bre{categorifiable} Let $R$ be a $\Zset_+$-ring categorified by an abelian 2-category $\mc{T}$. 
By \thref{fin-ab}(3) and \leref{prime-to-Serreprime}, the prime ideals of $R$ that are categorifiable are precisely the precisely the ones that are thick; that is the set 
\[
\Spec (R) \cap \mathrm{Serre}\mbox{-}\mathrm{Spec} (R). 
\]
\ere
\sectionnew{The Primitive Spectrum}
\label{7}
In this section we describe the relationship between the annihilation ideals of simple 2-representations of abelian 2-categories and 
the Serre prime ideals of these categories.
\subsection{2-representations}
\label{2reps}
Following Mazorchuk--Miemietz \cite{Mazorchuk2}, define a {\em{2-representation of a 2-category}} $\mc{T}$ to be a strict 
2-functor $\mc{F}$ from $\mc{T}$ to $\Cat$, the 2-category of all small categories. That is, $\mc{F}$ sends objects of $\mc{T}$ to small categories, 
1-morphisms of $\mc{T}$ to functors between categories, and 2-morphisms of $\mc{T}$ to natural transformations between functors.

Recall that the category of additive functors between two abelian categories has a canonical structure of an abelian category. 

\bde{exact-rep} A 2-representation $\mc{F}$ of a 2-category $\mc{T}$ will be called {\em{exact}}, if 
\begin{enumerate}
\item  $\mc{F}(A)$ is an abelian category for every object $A$ in $\mc{T}$;
\item $\mc{F}(f)$ is an additive functor for all 1-morphisms $f$ in $\mc{T}$; 
\item For any exact sequence of 1-morphisms in $\mc{T}$, 
\[
0 \to f \to g \to h \to 0,
\]
the sequence 
\[
0 \to \mc{F} (f) \to \mc{F} (g) \to \mc{F}(h) \to 0
\]
is an exact sequence of 1-morphisms in $\Cat$. 
\end{enumerate} 
\ede
Following, Mazorchuk, Miemietz, and Zhang \cite[Section 3.3]{Mazorchuk3}, we call a 2-representation $\mc{F}$ {\em{simple}} 
if the collection of categories 
\[
\{ \mc{F}(A) \mid A \in \mc{T}\}
\]
has no nonzero proper $\mc{T}$-invariant ideals. Such an ideal $X$ is a subset of the disjoint union of the set of morphisms of the categories  $\mc{F}(A)$ for $A \in \mc{T}$
with the following properties:
\begin{enumerate}
\item $ab$ and $ba$ are in $X$ for all $a \in X$ and all morphisms $b$ in $\mc{F}(A)$ such that the composition is well-defined;
\item $\mc{F}(f)(a) \in X$ for all $f \in \mc{T}_1$ and $a \in X$;
\item There is some morphism $a \in X$ which is not a zero morphism.
\end{enumerate}
\subsection{Annihilation ideals of 2-representations}
\label{annih-2r}
\bde{Annih}
Given an exact 2-representation $\mc{F}$ of the abelian 2-category $\mc{T}$, define its {\em{annihilation ideal}} $\Ann(\mc{F})$ to be the weak subcategory
of $\mc{T}$ having the same set of objects, set of 1-morphisms given by 
\[
\Ann(\mc{F})_1:=\{ f \in \mc{T}_1 \mid \mc{F}(f) \textrm{ is a zero functor}\},
\]
and set of 2 morphisms 
\[
\Ann(\mc{F})(f,g) := \mc{F}(f,g) \quad \mbox{for all} \quad f, g \in \Ann(\mc{F})_1.
\]
\ede
\ble{Annih} The annihilation ideal $\Ann(\mc{F})$ of every exact 2-representation $\mc{F}$ of an abelian 2-category $\mc{T}$ is a Serre 
ideal of $\mc{T}$.
\ele
\begin{proof} The proof is a direct verification of the necessary properties. 

To verify the ideal property of $\Ann(\mc{F})$, chose $f \in \Ann(\mc{F})_1$ and $g \in \mc{T}_1$ such that the composition is defined.
Then $\mc{F}(f)$ is a zero functor, and therefore, $\mc{F}(fg)=\mc{F}(f) \mc{F}(g)$ is also a zero functor. So, $fg \in \Ann(\mc{F})_1$. 
Likewise, $\mc{T}_1 \circ \Ann(\mc{F})_1 \subseteq \Ann(\mc{F})_1$.

To verify that $\Ann (\mc{F})(A_1, A_2)$ is a Serre subcategory of the abelian category $\mc{T}(A_1, A_2)$ for all object $A_1$ and $A_2$ of $\mc{T}$, 
consider an exact sequence $0 \to f \to g \to h \to 0$
in $\mc{T}(A_1, A_2)$. By \deref{exact-rep}(3), $0 \to \mc{F}(f) \to \mc{F}(g) \to \mc{F}(h) \to 0$ is an exact sequence in the abelian category 
of additive functors between the abelian categories $\mc{F}(A_1)$ and $\mc{F}(A_2)$. 

If $f, h \in \Ann(\mc{F})$, then $\mc{F}(f)$ and $\mc{F}(h)$ are both the zero functor and $\mc{F}(g)$ must also be the zero functor. Hence, $g \in \Ann (\mc{F})_1$. Likewise, assuming 
instead that $g \in \Ann(\mc{F})_1$, we get that $f, h \in \Ann(\mc{F})_1$. Hence, $\Ann(\mc{F})$ is a Serre ideal of $\mc{T}$. 
\end{proof}

Finally we have the following theorem, analogous to the relationship between prime ideals of rings and annihilators of simple representations, see 
e.g. \cite[Proposition 3.12]{Goodearl1}.

\bth{primitive}
Suppose $\mc{F}$ is a simple exact 2-representation of an abelian 2-category $\mc{T}.$ Then $\Ann(\mc{F})$ is a Serre prime ideal of $\mc{T}$. 
\eth
\begin{proof} We use the assumption of the simplicity of $\mc{F}$ to show that $\Ann(\mc{F})$ satisfies the condition in \thref{prime2}, form which we obtain that 
$\Ann(\mc{F})$ is a prime ideal of $\mc{T}$. The fact that $\Ann(\mc{F})$ is a Serre ideal was established in \leref{Annih}.

Suppose that $\mc{I}$ and $\mc{J}$ are thick ideals such that 
\[
\mc{I}_1 \circ \mc{J}_1 \subseteq \Ann(\mc{F})_1,
\] 
and neither $\mc{I}$ nor $\mc{J}$ is contained in $\Ann(\mc{F})$. 
Then we claim that the set
\[
X:=\{a  (\mc{F}(j) ( b)) c \mid a, b, c \textrm{ morphisms such that the composition is defined}, j \in \mc{J}_1 \}
\]
forms a nonempty $\mc{T}$-invariant ideal, contradicting the simplicity of $\mc{F}$.
It is clear that this set is an ideal, i.e., closed under composition on the left and right by any morphisms of $\Cat$ with appropriate source and target. 
We must show that it is invariant under $\mc{T}$, that it is nonzero, and that it is a proper subset of all morphisms of the categories $\mc{F}(A)$ 
for all objects $A$ of $\mc{T}$. 

First, assume that $g \in \mc{T}_1$. Then 
\[
\mc{F}(g) ( a (\mc{F}(j) (b)) c)= \mc{F}(g)(a) \mc{F}(g)(\mc{F}(j)(b)) \mc{F}(g)(c)=\mc{F}(g)(a) \mc{F}(gj)(b) \mc{F}(g)(c),
\]
which is clearly in $X$ whenever the composition is defined, since $gj \in \mc{J}_1$. Hence, $X$ is $\mc{T}$-invariant. 

Next, we show that $X$ is a proper subset. For all $i \in \mc{I}_1$ and $a (\mc{F}(j)(b))c \in X$, we have
\[
\mc{F}(i) (a (\mc{F}(j)(b))c)= \mc{F}(i)(a) \mc{F}(ij)(b) \mc{F}(i)(c)=\mc{F}(i)(a) 0 \mc{F}(i)(c)=0
\]
whenever the composition is defined. If $X$ equals the set of all morphisms of the collection of abelian categories $\{\mc{F}(A) \mid A \in \mc{T}\}$, then this would imply that 
$\mc{F}(i)$ is a zero functor for all $i \in \mc{I}_1$. Therefore, $\mc{I}_1 \subseteq \Ann(\mc{F})_1$. Applying \reref{incl-thick} and the assumption that $\mc{I}$ is a thick ideal 
gives that $\mc{I}$ is contained in $\Ann(\mc{F})$, which is a contradiction.
 
By a similar argument, one shows that $X$ contains nonzero morphisms. 
Since $\mc{J}$ is not contained in the annihilator of $\mc{F}$ by assumption, there is some $j \in \mc{J}_1$ such that $\mc{F}(j)$ is not the zero functor, 
and hence there is some morphism $b$ such that $\mc{F}(j)(b)$ is a nonzero morphism. 
Then by letting $a$ and $b$ be the appropriate identity morphisms, we see that $\mc{F}(j)(b)$ is a nonzero morphism in $X$.

Therefore, $X$ is a nonzero, proper $\mc{T}$-invariant ideal, which contradicts our assumption that $\mc{F}$ is simple. This gives that $\Ann(\mc{F})$ is a Serre prime ideal of $\mc{T}$.
\end{proof}

\bde{primitive}
The {\em{primitive spectrum}} of an abelian 2-category $\mc{T}$, denoted $\Prim (\mc{T})$, is the subset of $\mathrm{Serre}\mbox{-}\mathrm{Spec} (\mc{T})$ 
consisting of all primes $\mc{P}$ for which there exists a simple exact 2-representation $\mc{F}$ of $\mc{T}$ with $\mc{P}= \Ann (\mc{F})$. 
\ede
\sectionnew{Quantum Schubert cell algebras, canonical bases and prime ideals}
\label{8}
This section contains background material on quantum groups and quantum Schubert cell algebras, and their canonical bases defined by Kashiwara and Lusztig.
We recall facts about the 
homogeneous completely prime ideals of the quantum Schubert cell algebras and their relations to quantizations of Richardson varieties.
\subsection{Quantum groups, canonical bases and quantum Schubert cell algebras}
Let $\g$ be a (complex) symmetrizable Kac--Moody algebra with Cartan matrix $(a_{ij})_{i,j=1}^r$ and Cartan subalgebra 
$\t \subset \g$.
Denote the Weyl group of $\g$ by $W$. Let $\{\al_i \mid 1 \leq i \leq r \} \subset \t^*$ and $\{s_i \mid 1 \leq i \leq r \}$ be the sets of 
simple roots of $\g$ and simple reflections of $W$, respectively. Denote by $\{ \al_i \spcheck  \mid 1 \leq i \leq r \} \subset \t$ and  
$\{\vpi_i \mid 1 \leq i \leq r\} \subset \t^*$ the sets of simple coroots and fundamental weights of $\g$. Thus, 
$\langle \al_i\spcheck, \al_j \rangle = a_{ij}$. Let $(.,.)$ be the standard nondegenerate symmetric bilinear form on $\t^*$
satisfying 
\[
\langle \al_i \spcheck, \la \rangle =  \frac{ 2 ( \al_i,  \la)}{ ( \al_i , \al_i )} \; \; \mbox{for $\la \in \t^*$}
\quad \mbox{and} \quad \mbox{$( \al_i, \al_i ) = 2$ for short roots $\al_i$}. 
\]
Then
\[
d_i := \frac{(\al_i, \al_i)}{2} \in \Zset_+.
\]
Let 
\[
Q, \; \; P, \; \; P_+ \subset \t^*
\] 
be the root and weight lattices of $\g$, and the set of its dominant integral weights. Denote 
\[
P\spcheck := \{ h \in \t \mid \langle h, P \rangle \subset \Zset \} \subset \t \quad \mbox{and} \quad
Q_+:= \bigoplus \Zset_+ \al_i \subset \t^*.
\]

Let $U_q(\g)$ be the quantized universal enveloping algebra of $\g$ over $\Qset(q)$
with generators $e_i, f_i, q^h$ for $ 1 \leq i \leq r$, $h \in P\spcheck$ and relations as in \cite{Kashiwara1}.
We will use the Hopf algebra structure of $U_q(\g)$ with coproduct given by
\begin{equation}
\label{coproduct}
\Delta(e_i) = e_i \otimes 1 + q^{d_i \al_i\spcheck} \otimes e_i, \quad
\Delta(f_i) = f_i \otimes q^{- d_i \al_i\spcheck} + 1 \otimes f_i, \quad
\Delta(q^h) = q^h \otimes q^h
\end{equation}
for $h \in P\spcheck$, $1 \leq i \leq r$. Let  $U^\pm_q(\g)$ and $U^0_q(\g)$ be the unital subalgebras of $U_q(\g)$ 
generated by $\{e_i \mid 1 \leq i \leq r \}$ (resp. $\{f_i \mid 1 \leq i \leq r \}$) and $\{ q^h \mid h \in P\spcheck \}$. 
Denote the (symmetric) $q$-integers and factorials 
\[
q_i := q^{d_i}, \quad [n]_i := \frac{q_i^n - q_i^{-n}}{q_i - q_i^{-1}} \quad \mbox{and} \quad
[n]_i! := [1]_i \ldots [k]_i. 
\]

Denote by $*$ and $\varphi$ the $\Qset(q)$-linear anti-automorphisms of $U_q(\g)$ defined by
\begin{align*}
&e_i^* := e_i, \quad f_i^* := f_i, \quad (q^h)^* := q^{-h}, \quad \mbox{and} \\
&\varphi(e_i) := f_i, \quad \varphi(f_i) := e_i, \quad \varphi(q^h) := q^h
\end{align*}
for $1 \leq i \leq r$, $h \in P\spcheck$. 
The composition $\varphi^* := \varphi \circ * = * \circ \varphi$, which is a $\Qset(q)$-linear automorphism of $U_q(\g)$, 
satisfies
\[
\varphi^*(e_i) = f_i, \quad \varphi^*(f_i) = e_i \quad \mbox{and} \quad \varphi^*(q^h) = q^{-h}.
\]
This composition is denoted by ${}\spcheck$ in \cite{Kashiwara5}; we will use the above notation to avoid interference with later notation.

The Hopf algebra $U_q(\g)$ is graded by the root lattice $Q$ by setting 
\begin{equation}
\label{grading}
\deg e_i =  \al_i, \quad \deg f_i = - \al_i, \quad 
\deg q^h =0.
\end{equation}
The homogeneous components of a subspace $Y$ of $U_q(\g)$ of degree $\gamma \in Q$
will be denoted by $Y_\gamma$.
Denote by $e''_i$ the $\Qset(q)$-linear skew-derivations of $U_q^-(\g)$ such that
\[
e''_i (f_j) = \delta_{ij} \quad \mbox{and} \quad
e''_i(x y) = e''_i(x) y + q^{ - ( \al_i, \gamma )} x e''_i(y) \quad
\mbox{for} \quad x \in U_q^-(\g)_\gamma, y \in U_q^-( \g).
\]
Kashiwara's (nondegenerate, symmetric) bilinear form $(-, -)_K : U_q^-(\g) \times U_q^-(\g) \to \Qset(q)$ is defined by
\[
(1, 1)_K =1 \quad \mbox{and} \quad (f_i x, y)_K = (x, e''_i(y))_K 
\]
for all $1 \leq i \leq r$ and $x, y \in U_q^-(\g)$. 

\bre{Kashiwara-form} This differs slightly from the conventional choice for Kashiwara's form $\langle -, - \rangle_K : U_q^-(\g) \times U_q^-(\g) \to \Qset(q)$, 
which is defined by
\[
\langle 1, 1 \rangle_K =1 \quad \mbox{and} \quad \langle f_i x, y \rangle_K = \langle x, e'_i(y) \rangle_K 
\]
for all $1 \leq i \leq r$ and $x, y \in U_q^-(\g)$ in terms of the  
$\Qset(q)$-linear skew-derivations $e'_i$ of $U_q^-(\g)$
given by
\[
e'_i (f_j) = \delta_{ij} \quad \mbox{and} \quad
e'_i(x y) = e'_i(x) y + q^{( \al_i, \gamma )} x e'_i(y) \quad
\mbox{for} \quad x \in U_q^-(\g)_\gamma, y \in U_q^-( \g).
\]
The two forms are related by
\begin{equation}
\label{bar}
(x,y)_K = \overline{\langle \overline{x}, \overline{y} \rangle}_K \quad \mbox{for} \quad x, y \in U_q^-(\g),
\end{equation}
where $x \mapsto \overline{x}$ denotes the $\Qset$-linear automorphism of $\Qset(q)$ given by $\bar{q} = q^{-1}$ and 
the bar involution of $U_q(\g)$ (the skew-linear automorphism of $U_q(\g)$ given by $\overline{f}_i = f_i$).
Using \eqref{bar}, one converts dualization results with respect to one form to such results for the other.
\ere

Let $\AA := \Zset[q^{\pm 1}]$ and $U_\AA^\pm(\g)$ be the (divided power) {\em{integral forms}} of $U_q^\pm(\g)$, which are the 
$\AA$-subalgebras of $U_q^\pm(\g)$ generated by $e_i^{(k)} = e_i^k/[k]_i!$
and $f_i^{(k)} = f_i^k/[k]_i!$ for $1 \leq i \leq r$, $k \in \Zset_+$,
respectively. The {\em{dual integral form}} $U_\AA^-(\g)\spcheck$ of $U_q^-(\g)$ 
is the $\AA$-subalgebra given by
\[
U_\AA^-(\mathfrak{g})\spcheck=\{x \in U_q^-(\mathfrak{g}) \mid ( x, U_\AA^-(\mathfrak{g}))_K \subset \AA\}.
\]
 
Kashiwara \cite{Kashiwara1} and Lusztig \cite{Lusztig1} defined the {\em{canonical/lower global basis}} of $U_\AA^\pm (\g)$ and the {\em{dual canonical/upper global basis}} of 
$U_\AA^\pm(\g)\spcheck$. These bases have a number of remarkable properties; for instance, they descend to bases of integrable highest weight modules 
by acting on highest weight vectors. We will denote by ${\BB}_\pm^{\low}$ the lower global basis of $U_\AA^\pm(\g)$ and by ${\BB}_-^{\up}$ the 
upper global basis of $U_\AA^-(\g)$.  

The lower global basis ${\BB}_-^{\low}$ and the upper global basis ${\BB}_-^{\up}$ form a pair of dual bases of $U_\AA^-(\g)$ and 
$U_\AA^-(\g)\spcheck$ with respect to the pairing $(-,-)_K$. For $b \in {\BB}_-^{\low}$, denote by 
$b\spcheck \in {\BB}_-^{\up}$ the corresponding dual element, so 
\begin{equation}
\label{dualb}
( b, c\spcheck )_K = \delta_{b,c} \quad \mbox{for} \quad
b, c \in {\BB}_-^{\low}.
\end{equation}

The lower global bases ${\BB}_\pm^{\low}$ satisfy the invariance properties
\begin{equation}
\label{invarian}
({\BB}_\pm^{\low})^* = {\BB}_\pm^{\low} \quad \mbox{and} \quad \varphi({\BB}_\pm^{\low}) = \varphi^*({\BB}_\pm^{\low}) = {\BB}_\mp^{\low},
\end{equation} 
see \cite[Theorem 2.1.1]{Kashiwara4}, \cite[Theorem 4.3.2]{Kashiwara5} and \cite[Theorem 8.3.4]{KKKO}.

To each Weyl group element $w$, one associates the quantum Schubert cell algebras $U^-_q[w] \subseteq U_q^-(\g)$. They can be defined in 
two ways. Starting from a reduced expression 
\[
w = s_{i_1} \ldots s_{i_N}
\]
of $w$, consider the roots
\[
\be_1 := \al_{i_1}, \be_2 := s_{i_1} (\al_{i_2}), \ldots, 
\be_N := s_{i_1} \ldots s_{i_{N-1}}(\al_{i_N})
\]
and the root vectors
\begin{equation}
\label{roots}
\{ f_{\be_j} := 
T_{i_1} \ldots T_{i_{j-1}} 
(f_j) \mid 
1 \leq j \leq N \}
\end{equation}
using Lustig's braid group action \cite{Lusztig2,Jantzen1} on $U_q(\g)$. 
De Concini, Kac, and Procesi \cite{DKP}, and Lusztig \cite[\S 40.2]{Lusztig2} defined the 
algebra $U^-_q[w]$ as the unital $\Qset(q)$-subalgebra of $U_q^-(\g)$ with generating set \eqref{roots}, 
and proved that this is independent on the choice of reduced expression of $w$. 
Berenstein and Greenstein \cite{BerGreen} conjectured that
\[
U^-_q[w] = U^-_q(\g) \cap T_w ( U^+_q(\g)),
\]
and Kimura \cite{Kimura2} and Tanisaki \cite{Tani} proved this property. It can be used as a second definition of the algebras $U^-_q[w]$. 
Kimura proved \cite[Theorem 4.5]{Kimura1} that
\begin{equation}
\label{uppergbBw}
{\BB}_-^{\up}[w] := {\BB}_-^{\up} \cap U^-_q[w] 
\end{equation}
is an $\AA$-basis of the $\AA$-algebra
\[
U^-_\AA[w]\spcheck := U^-_q[w] \cap U^-_\AA(\g)\spcheck. 
\]
We will refer to this algebra as to the {\em{dual integral form}} of $U^-_q[w]$. 
The set ${\BB}_-^{\up}[w]$ is called the {\em{upper global basis}} of $U^-_\AA[w]\spcheck$.
\subsection{Homogeneous completely prime ideals of the algebras $U^-_q[w]$}
Denote the Hopf subalgebras $U^{\geq 0}:= U_q^+(\g) U_q^0(\g)$ and $U^{\leq 0} :=  U_q^-(\g) U_q^0(\g)$ of $U_q(\g)$. 
The Rosso-Tanisaki form
\[
( -, -)_{RT} : U^{\leq 0} \times U^{\geq 0} \to \Qset(q^{1/d})
\]
(for an appropriate $d \in \Zset_+$) is the Hopf algebra pairing satisfying
\[
(y, x x')_{RT} = (\Delta(y), x' \otimes x)_{RT}, \quad (y y' , x)_{RT} = (y \otimes y', \Delta(x))_{RT}
\]
for $y, y' \in U^{\leq 0}$, $x, x' \in U^{\geq 0}$, and normalized by
\[
(f_i, e_j)_{RT} = \delta_{ij}, \quad (q^h, q^{h'})_{RT} = q^{- ( h, h' )}, \quad
(f_i, q^{h'})_{RT}= (q^h, e_i)_{RT}=0
\]
for $1 \leq i,j \leq r$, $h,h' \in P\spcheck$. This is a slightly different normalization than the usual one \cite[Eq. (6.12)(2)]{Jantzen1}
needed in order to match this form to Kashiwara's one. The two normalizations for $(-,-)_{RT}$ are related to each other 
by a Hopf algebra automorphism of $U^{\leq 0}$ coming from the torus action associated to its $Q$-grading.

For $\gamma \in Q_+$ let
\[
\{ x_{\gamma, i} \} \quad \mbox{and} \quad \{ y_{\gamma, i} \}
\] 
be a set of dual bases of $(U^-_q(\g))_{-\gamma}$ and $(U^+_q(\g))_{\gamma}$ with respect to $(-,-)_{RT}.$ The {\em{quasi-$R$-matrix}} of $U_q(\g)$
is 
\[
\mc{R} := \sum_{\gamma \in Q_+} \sum _i y_{\gamma, i} \otimes x_{\gamma, i} \in U_q^+(\g) \widehat{\otimes} U_q^-(\g)
\]
where the completed tensor product is with respect to the descending filtration \cite[\S 4.1.1]{Lusztig2}.

For $\la \in P_+$, we will denote by $V(\la)$ the irreducible $U_q(\g)$-module of highest weight $\la$ and by 
$V(\la)^\circ$ its restricted dual
\begin{equation}
\label{dualV}
V(\la)^\circ := \oplus_{\nu \in P} (V(\la)_\nu)^*
\end{equation}
where 
\[
V(\la)_\nu := \{ v \in V(\la) \mid q^h v = q^{ \langle \nu, h \rangle} v, \; \forall h \in P\spcheck \} \quad \mbox{for} \quad
\nu \in P
\]
are the (finite dimensional) weight spaces of $V(\la)$. Let $v_\la$ be a fixed highest weight vector of $V(\la)$. 
Denote by ${\BB}(\la)^{\low}$ the {\em{lower global basis}} of (the integral form $U_\AA^-(\g) v_\la$) of $V(\la)$.
It is an $\AA$-basis of $U_\AA^-(\g) v_\la$ and a $\Qset(q)$-basis of $V(\la)$.
For $w \in W$, let $v_{w \la}$ be the unique element of $V(\la)_{w \la}$ which belongs to ${\BB}(\la)^{\low}$. 
Let 
\[
V_w^\pm(\la) := U^\pm_q(\g) v_{w \la} \subseteq V(\la)
\]
be the associated {\em{Demazure modules}}. 
For $v \in V(\la)$ and $\xi \in V(\la)^*$, denote the 
corresponding {\em{matrix coefficient}} of $V(\la)$ considered as a functional on $U_q(\g)$:
\[
c_{\xi, v} \in (U_q(\g))^* \quad \mbox{given by} 
\quad
c_{\xi, v}( x) = \xi( x \cdot v) \quad \mbox{for} \; \; x \in U_q(\g).
\]

A subspace $U$ of $U_q(\g)$ will be called {\em{homogeneous}} if 
\[
U = \bigoplus_{\gamma \in Q} U_\gamma \quad \mbox{where} \quad U_\gamma := U \cap U_q(\g)_\gamma.
\]  
\bth{Uw-id} \cite[Theorem 3.1(a)]{Y-plms} Let $\g$ be a symmetrizable Kac--Moody algebra and $w \in W$ be a Weyl group element. For all $u \in W$ such that $u \leq w$, 
the set 
\[
I_w(u) = \{ \langle c_{\xi, v_{w \la}} \otimes \id, \mc{R} \rangle^* \mid \xi \in V(\la)^\circ, 
\xi \perp V_u^-(\la), \la \in P_+\}
\]
is a homogeneous completely prime ideal of $U^-_q[w]$.
\eth

The original proof of this theorem that was given in \cite{Y-plms} extensively used the works of Joseph \cite{Joseph-JA,Joseph-book} and Gorelik \cite{Gorelik}
which were written for the case when $\g$ is finite dimensional (complex simple Lie algebra). In \cite[Lemmas 3.4 and 3.5]{Y-GMJ} 
a simplified, self-contained proof was given that established the 
validity of the theorem when $U_q(\g)$ is defined over a base field $\kk$ of arbitrary characteristic; 
the needed results of Joseph \cite{Joseph-JA,Joseph-book} and Gorelik \cite{Gorelik} were reproved in this setting. 
The proofs of these facts in \cite{Y-GMJ} work in the full 
generality of symmetrizable Kac--Moody algebras. 

The Rosso-Tanisaki form $(-,-)_{RT}$ satisfies
\[
( y^*, x^* )_{RT} = (y,x)_{RT} \quad \mbox{for} \quad y \in U^{\leq 0}, x \in U^{\geq 0}, 
\]
\cite[Lemma 6.16]{Jantzen1}. Therefore, 
\begin{equation}
\label{Rstst}
\mc{R}^{*\otimes *} =  \mc{R}, 
\end{equation}
and thus, the ideals $I_w(u)$ are also given by 
\[
I_w(u) = \{ \langle (c_{\xi, v_{w \la}} \circ *) \otimes \id, \mc{R} \rangle \mid \xi \in V(\la)^\circ, 
\xi \perp V_u^-(\la), \la \in P_+\}.
\]

We will need the following relation between the bilinear forms $(-,-)_K$ and $(-, -)_{RT}$, and the corresponding expression for $\mc{R}$ in terms of global bases. 
\bpr{forms} Let $\g$ be a symmetrizable Kac--Moody algebra. For all $x_1, x_2 \in U_q^-(\g)$, we have
\begin{equation}
\label{relat-bilinf}
(x_1, \varphi^*(x_2))_{RT} = (x_1, x_2)_K.
\end{equation}
The quasi-$R$-matrix of $U_q(\g)$ is given by
\begin{equation}
\label{Rmatr}
\mc{R} = \sum_{b \in {\BB}_-^{\low}} \varphi^*(b) \otimes b\spcheck
= \sum_{b \in {\BB}_-^{\low}} \varphi(b) \otimes (b\spcheck)^*,
\end{equation}
recall \eqref{dualb}.
\epr
\begin{proof} The $Q$-grading of $U^{\leq 0}$ specializes to a $\Zset_+$-grading via the group homomorphism $Q \to \Zset$ given by $\al_i \mapsto -1$.
The corresponding graded components will be denoted by $(U^{\leq 0})_l$. Set
\[
(U^{\leq 0})_{\geq l} := (U^{\leq 0})_l \oplus (U^{\leq 0})_{l+1} \oplus \ldots
\] 
For $x := f_{i_1} \ldots f_{i_k}$ and $h := d_{i_1} \al_{i_1}\spcheck + \cdots + d_{i_k} \al_{i_k}\spcheck$, we have
\[
\Delta(x) - x \otimes q^{-h} - \sum_{j=1}^k q^{ ( \al_{i_j}, \al_{i_1} + \ldots + \al_{i_{j-1}} ) }
f_{i_1} \ldots f_{i_{j-1}} f_{i_{j+1}} \ldots f_{i_k} \otimes f_{i_j} q^{-h + d_{i_j} \al_{i_j}\spcheck} \in 
U^{\leq 0} \otimes (U^{\leq 0})_{\geq 2},
\]
i.e.,
\[
\Delta(x) - x \otimes q^{-h} - \sum_{i=1}^r e''_i(x) \otimes f_i q^{- h + d_i \al_i\spcheck} \in U^{\leq 0} \otimes (U^{\leq 0}){\geq 2}.
\]
This property, and the two properties of the Rosso-Tanisaki form
\begin{align*}
&( (U^{\leq 0})_\gamma, (U^{\geq 0})_\nu )_{RT} =0 \quad \mbox{for} \quad \gamma + \nu \neq 0, 
\\
&(x q^h, y q^{h'})_{RT} = q^{ - ( h, h' )} (x, y)_{RT} \quad \mbox{for} \quad x \in U_q^-(\g), y \in U_q^+(\g), h, h' \in P\spcheck 
\end{align*}
(see \cite[Eqs. 6.13(1)-(2)]{Jantzen1}) imply that the bilinear form $ \langle -, - \rangle$ on $U_q^-$ 
given by
\[
\langle x , y \rangle := ( x, \varphi^*(y))_{RT} 
\]
satisfies $\langle x,  f_i y \rangle = \langle e''_i(x), y \rangle$ for $x, y \in U_q^-(\g)$ and $1 \leq i \leq r$. The uniqueness 
property of the Kashiwara form implies that this form equals $( - , - )_K$, which proves \eqref{relat-bilinf}. 

The invariance property \eqref{invarian}, the relation \eqref{relat-bilinf} between the bilinear forms $(-,-)_{RT}$ and $(-,-)_K$, 
and the orthogonality property \eqref{dualb} imply that
\[
\{ \varphi^*(b) \mid b \in {\BB}_-^{\low} \} \quad \mbox{and} \quad
\{ b\spcheck \mid b \in {\BB}_-^{\low} \}
\]
are a pair of dual bases of $U_q^+(\g)$ and $U_q^-(\g)$ with respect to the pairing $(-,-)_{RT}$. This gives the first equality in
\eqref{Rmatr}. The second equality follows from \eqref{Rstst} and the first invariance property in \eqref{invarian}.
\end{proof}
\subsection{Quantizations of Richardson varieties}
Let $G$ be the Kac--Moody group (over $\Cset$) corresponding to $\g$.  Let $B_{\pm}$ be opposite Borel subgroups of $G$. 
For $u, w \in W$, the {\em{open Richardson variety}} associated to the pair $(u,w)$ is the locally closed subset of the flag variety $G/B_+$ defined by 
\[
R_{u,w}:=(B_- u B_+)/B_+ \cap (B_+ w B_+)/B_+.
\]
It is nonempty if and only if $u \leq w$ in which case it has dimension $\ell(w) - \ell(u)$ (in terms of the standard length function $\ell : W \to \Zset_+$).
We have the stratifications of $G/B_+$ into unions of Schubert cells
\[
G/B_+=\coprod_{w \in W} (B_+ w B_+)/B_+= \coprod_{u \in W} (B_- u B_+)/B_+
\]
and open Richardson varieties
\[
G/B_+=\coprod_{\substack{ u \leq w\\ u,w \in W}} R_{u,w}.
\]
Denote the closure
\[
\overline{R}_{u,w} := \mathrm{Cl}_{(B_+ w B_+)/B_+} (R_{u,w})
\]
of $R_{u,w}$ in the Schubert cell $(B_+ w B_+) /B_+$. 

For $w \in W$, define the dual extremal vectors
\[
\xi_{w \la} \in V(\la)_{w \la}^* \quad \mbox{by} \quad \xi_{w \la} (v_{w \la}) =1
\]
(keeping in mind that $\dim V(\la)_{w \la} =1$). 
Denote the image of the corresponding extremal matrix coefficient in $U^-_q[w]$:
\[
\Delta_{\la, w \la}:= \langle c_{\xi, v_{w \la}} \otimes \id, \mc{R} \rangle^* \in U^-_q[w] \quad \mbox{for} \quad \la \in P_+.
\]
\bpr{quant}
For all symmetrizable Kac--Moody algebras $\g$ and $u \leq w \in W$, the factor ring $U^-_q[w]/I_w(u)$ 
is a quantization of the coordinate ring $\Cset[\overline{R}_{u,w}]$ of the closure of the open Richardson variety $R_{u,w}$ in the Schubert cell $(B_+ w B_+)/B_+$. The localization 
\[
(U^+[w]/I_w(u))[\Delta_{\vpi_i, w \vpi_i}^{-1}, 1 \leq i \leq r]
\]
of this ring is a quantization of the coordinate ring $\Cset[R_{u,w}]$. 
\epr
These facts were stated in \cite[pp. 274-275]{Y-contmath} for finite dimensional complex simple Lie algebras $\g$, but the proofs given there 
carry over to the symmetrizable Kac--Moody case directly.
\sectionnew{Categorifying Richardson Varieties}
\label{9}
In this section, we prove that the ideals $I_w(u) \cap U^-_\AA[w]\spcheck$ of $U^-_\AA[w]\spcheck$ are Serre completely prime ideals for all symmetric Kac--Moody algebras $\g$
and $u \leq w \in W$.
We then use \thref{fin-ab}(4) to construct a (domain) multiring category which categorifies the quantization of the coordinate ring of the
closure of the open Richardson variety $R_{u,w}$ in the Schubert cell $(B_+ w B_+)/B_+$. This category is obtained as a factor of a multiring category 
consisting of graded, finite dimensional representations of the corresponding KLR algebras.
\subsection{The categorifications of $U_\AA^-(\mathfrak{g})\spcheck$ and $U_\AA^-[w]\spcheck$, and relations to dual canonical bases}
For each symmetrizable Kac--Moody algebra $\g$, Khovanov, Lauda  \cite{Khovanov1} and  Rouquier \cite{Rouquier0} 
defined a family of (graded) quiver Hecke algebras over a base field $\kk$,
which we will call KLR algebras. They proved that the category $\mc{C}$ which is the direct sum of the categories of finite dimensional graded modules of 
the KLR algebras associated to $\g$ has the following properties:
\bth{KLR} {\em{(}}\hspace{-0.04cm}Khovanov--Lauda \cite{Khovanov1} and  Rouquier \cite{Rouquier0}{\em{)}} 
For each symmetrizable Kac--Moody algebra $\g$ and base field $\kk$, $\mc{C}$ is a 
$\kk$-linear multiring category such that
\begin{equation}
\label{K0KLR}
K_0(\mc{C}) \cong U_\AA^-(\g)\spcheck.
\end{equation}
The action of $q$ on the right hand side and the shift of grading autoequivalence of $\mc{C}$ are related via
\[
[M(k)] = q^k [M] \quad \mbox{for all objects $M$ of $\mc{C}$}.
\]
\eth 
The theorem implies that for every symmetrizable Kac--Moody algebra $\g$, $U_\AA^-(\g)\spcheck$ is a $\Zset_+$-ring
with positive $\Zset_+$-basis $\{[M]\}$ where $M$ runs over the isomorphism classes of the simple objects of $\mc{C}$. 
(Here we disregard the structure of $U_\AA^-(\g)\spcheck$ as an $\AA$-algebra and view it just as a ring.) 
For symmetric Kac--Moody algebras $\g$, the relation between this basis and the upper global basis of $U_\AA^-(\g)\spcheck$
is given by the next theorem. For it we recall that the dual of each graded finite dimensional representation of a KLR algebra 
has a canonical structure of a KLR module which is also graded, finite dimensional. This gives a canonical duality endofunctor of $\mc{C}$. 

\bth{KLR-cano} {\em{(}}\hspace{-0.05cm}Varagnolo--Vasserot \cite{Varagnolo1} and Rouquier \cite{Rouquier1}{\em{)}} 
For each symmetric Kac--Moody algebra $\g$ and base field $\kk$ of characteristic 0, 
under the isomorphism \eqref{K0KLR}, the upper global basis corresponds to the set of isomorphism classes of the the self-dual simple modules in 
the category $\mc{C}$.
\eth
The theorem implies that in these cases $U_\AA^-(\g)\spcheck$ is a $\Zset_+$-ring with positive $\Zset_+$-basis 
\[
q^\Zset {\BB}_-^{\up}.
\]

For each symmetric Kac--Moody algebra $\g$ and $w \in W$, in \cite[\S 11.2]{KKKO} Kang, Kashiwara, Kim and Oh constructed a monoidal subcategory 
$\mc{C}_w$ of $\mc{C}$ as the smallest monoidal Serre subcategory closed under shifts, containing a certain set of simple self-dual modules of the KLR algebras
(\cite[Definition 11.2.1]{KKKO}) and, using \cite{GLS}, they proved:
\bth{catUw} {\em{(}}\hspace{-0.05cm}Kang--Kashiwara--Kim--Oh \cite{KKKO}{\em{)}} For each symmetric Kac--Moody algebra $\g$, 
\begin{equation}
\label{K0wiso}
K_0( \mc{C}_w) \cong U^-_\AA[w]\spcheck.
\end{equation}
\eth
Combining Theorems \ref{tKLR-cano} and \ref{tcatUw} gives that, under the isomorphism \eqref{K0wiso}, the elements of the upper global basis
${\BB}_-^{\up}[w]\spcheck$ of $U_\AA^-[w]\spcheck$ (recall \eqref{uppergbBw})
correspond to the isomorphism classes of the simple self-dual objects of $\mc{C}_w$. 
In particular, $U_\AA^-[w]\spcheck$ is a $\Zset_+$-ring with a positive $\Zset$-basis
\[
q^\Zset {\BB}_-^{\up}[w].
\] 
\subsection{Serre completely prime ideals of the $\Zset_+$-rings $U_\AA^-[w]$ and the multiring categories $\mc{C}_w$}
For $u \leq w$, denote the ideals
\[
I_w(u)_\AA\spcheck := I_w(u) \cap U_\AA^-[w]\spcheck
\]
of $U_\AA^-[w]$. \thref{Uw-id} implies that for every symmetrizable Kac--Moody algebra $\g$, $I_w(u)_\AA\spcheck$ are 
completely prime ideals of $U_\AA^-[w]$ in the classical sense. The following is the main result of this section.

\bth{cat-Rich} (1) Let $\g$ be a symmetrizable Kac--Moody algebra and $u \leq w \in W$.
The ideal $I_w(u)_\AA\spcheck$ has an $\AA$-basis given by
\[
{\BB}_-^{\up}[w] \cap I_w(u)_\AA\spcheck.
\]
Furthermore, it is a Serre completely prime ideal of the $\Zset_+$-ring $U_\AA^-[w]\spcheck$.

Denote by $X_w(u)$ the set of isomorphism classes of self-dual simple objects $M$ of $\mc{C}_w$ such that $[M] \in I_w(u)_\AA\spcheck$. Let 
\[
\mc{I}_w(u) := \mc{S}(X_w(u)[k], k \in \Zset)
\]
be the full subcategory of $\mc{C}_w$ whose objects have Jordan-H\"older series with simple subquotients isomorphic to shifts of objects in $X_w(u)$  
as in \leref{fin-ab}(2). 

(2) Let $\g$ be a symmetric Kac--Moody algebra and $u \leq w \in W$.
For all base fields $\kk$ of characteristic 0, $\mc{I}_w(u)$ are Serre completely prime ideals of the $\kk$-linear multiring category $\mc{C}_w$. For the corresponding 
Serre quotient $\mc{C}_w/ \mc{I}_w(u)$, we have 
\[
K_0( \mc{C}_w / \mc{I}_w(u) ) \cong U_\AA^-[w]\spcheck/I_w(u)_\AA\spcheck.
\]
\eth
By the first part of \thref{cat-Rich}(1), 
\[
\left( U_\AA^-[w]\spcheck/I_w(u)_\AA\spcheck \right) \otimes_\AA \Qset(q) \cong U^-[w]/I_w(u)
\]
and by \prref{quant}, $U^-[w]/I_w(u)$ is a quantization of the coordinate ring $\Cset[\overline{R}_{u,w}]$ of the closure 
of the Richardson variety $R_{u,w}$ in the Schubert cell $(B_+ w B_+)/B_+$. This fact and \thref{cat-Rich}(2) imply that the Serre quotient
$\mc{C}_w/\mc{I}_w(u)$, which is a domain in the sense of \deref{domain}, 
is a monoidal categorification of the quantization of $\Cset[\overline{R}_{u,w}]$.
\subsection{Proof of \thref{cat-Rich}}
Recall that ${\BB}(\la)^{\low}$ denoted the lower global basis of the irreducible module $V(\la)$ 
for $\la \in P_+$. We will need two facts about the lower global bases of Demazure modules proved by Kashiwara:
\bth{basis1} {\em{(}}\hspace{-0.05cm}Kashiwara \cite{Kashiwara4}{\em{)}}
For every symmetrizable Kac-Moody algebra $\g$ and dominant integral weight $\la \in P_+$, the intersection 
\[
{\BB}_w^\pm(\la)^{\low}:={\BB}(\la)^{\low} \cap V^\pm_w(\la)
\]
is a $\Qset(q)$-basis of the Demazure module $V^\pm_w(\la)$.
\eth

The sets ${\BB}_w^\pm(\la)^{\low}$ are called the lower global bases of the Demazure modules $V^\pm_w(\la)$.
The plus case was proved in \cite[Proposition 3.2.3(i)]{Kashiwara4} and the minus in \cite[Proposition 4.1]{Kashiwara4}. 
The following theorem describes the relationship between the canonical/lower global bases ${\BB}_w^+(\la)$ 
of the Demazure modules and the action of the canonical/lower global bases of $U_\AA^+(\g)$ acting on the 
corresponding extremal weight vectors.
\bth{basis2} {\em{(}}\hspace{-0.05cm}Kashiwara \cite{Kashiwara5,Kashiwara6}{\em{)}} Let $\g$ be a symmetrizable Kac-Moody algebra $\g$ 
and $\la \in P_+$ be a dominant integral weight. Denote the subset
\[
\mathfrak{B}_w^+(\la)^{\low}:=\{ b \in {\BB}_+^{\low} \mid b \cdot v_{w \la} \neq 0 \} 
\]
of the lower global basis of $U^+_\AA(\g)$. Then there is a bijection between this set and the lower global basis of the Demazure module 
$V_w^+(\la)$ given by
\[
\eta_w : \mathfrak{B}_w^+(\la)^{\low} \stackrel{\cong}\longrightarrow {\BB}_w^+(\la)^{\low} \quad \mbox{given by} \quad \eta_w(b) := b \cdot v_{w \la}.
\]
\eth
The corresponding fact to this theorem for the negative Demazure modules (where everywhere plus is replaced by minus) was proved in 
\cite[Proposition 4.1]{Kashiwara4}.

The following proposition is a stronger form of the statement of the first part of \thref{cat-Rich}(1).
\bpr{Iuw} For all symmetrizable Kac--Moody algebras $\g$, and $u \leq w \in W$, the ideal $I_{w}(u)$ of the quantum Schubert cell algebra 
$U^-_q[w]$ has a $\Qset(q)$-basis given by
\[
\bigcup_{\la \in P_+} \big\{ b\spcheck \mid b \in \varphi^{-1} \eta_w^{-1} \big( {\BB}^+_w(\la)^{\low} \backslash {\BB}^-_u(\la)^{\low} \big) \big\}.
\]
\epr
\begin{proof} For $\la \in P_+$, consider the basis of $V(\la)^\circ$ (cf. \eqref{dualV}) which is dual to the lower global basis $B(\la)^{\low}$ of $V(\la)$. 
Given $v \in B(\la)^{\low}$, denote by $v\spcheck$ the corresponding dual element, so
\[
v_1\spcheck (v_2) = \delta_{v_1 v_2} \quad \mbox{for} \quad v_1, v_2 \in B(\la)^{\low}.
\] 
 
\thref{basis1} implies that
\begin{multline*}
\{ \xi \in V(\la)^\circ \mid \xi \perp V_u^-(\la) \} = 
\Span_{\Qset(q)} \{ {\BB}^+_w(\la)^{\low} \backslash {\BB}^-_u(\la)^{\low} \}  \\
\oplus 
\Span_{\Qset(q)} \{ {\BB}(\la)^{\low} \backslash ({\BB}^+_w(\la)^{\low} \cup {\BB}^-_u(\la)^{\low} ) \}.
\end{multline*}
For $v \in {\BB}(\la)^{\low} \backslash ( {\BB}^+_w(\la)^{\low} \cup {\BB}^-_u(\la)^{\low})$, we have $v\spcheck \perp V_w^+(\la)$, and thus,
\[
\langle c_{v\spcheck, v_{w \la}} \otimes \id, \mc{R} \rangle = 0.
\] 
Therefore, the subspace
\[
\{ \langle c_{\xi, v_{w \la}} \otimes \id, \mc{R} \rangle^* \mid
\xi \in V(\la)^\circ \mid \xi \perp V_u^-(\la) \} \subset I_w(u)
\]
is spanned by
\[
\{ \langle c_{v\spcheck, v_{w \la}} \otimes \id, \mc{R} \rangle^* \mid
v \in {\BB}^+_w(\la)^{\low} \backslash {\BB}^-_u(\la)^{\low} \}.
\]
The proposition now follows from the identity
\begin{equation}
\label{c-to-B}
\langle c_{v\spcheck, v_{w \la}} \otimes \id, \mc{R} \rangle^* = \big( \varphi^{-1} \eta_w^{-1}(v) \big)\spcheck
\quad \mbox{for} \quad 
v \in {\BB}^+_w(\la)^{\low}.
\end{equation}
To show this, first note that 
\thref{basis2} implies that for $v \in {\BB}^+_w(\la)^{\low}$ and $b \in {\BB}_+^{\low}$,
\[
\langle v\spcheck, b \cdot v_{w \la} \rangle = 
\begin{cases}
1, & \mbox{if} \; \; b = \eta_w^{-1}(v) 
\\
0, & \mbox{otherwise}.
\end{cases} 
\]
Using this and the second part of \prref{forms}, for $v \in  {\BB}^+_w(\la)^{\low}$, we obtain
\[
\langle c_{v\spcheck, v_{w \la}} \otimes \id, \mc{R} \rangle^* =  \sum_{b \in {\BB}_-^{\low} } 
\langle v\spcheck, \varphi(b) \cdot v_{w \la} \rangle b\spcheck
= \big( \varphi^{-1} \eta_w^{-1}(v) \big)\spcheck,
\]
which shows \eqref{c-to-B} and completes the proof of the proposition.
\end{proof}

\begin{proof}[Proof of \thref{cat-Rich}] (1) The first statement in part (1) follows from \prref{Iuw}. 
By \thref{Uw-id}, $I_w(u)$ is a completely prime ideal of $U^-_q[w]$, and therefore, the contraction
$I_w(u) \cap U^-_\AA[w]$ is a completely prime ideal of $U^-_\AA[w]$. The ideal
$I_w(u) \cap U^-_\AA[w]$ has a $\Zset$-basis consisting of elements that belong to 
$q^\Zset {\BB}^{\up}_-[w]$, which, by \thref{KLR-cano}, is precisely the positive basis of the $\Zset_+$-ring $U^-_\AA[w]$.
Now we apply \leref{prime-to-Serreprime}, which gives that
$I_w(u) \cap U^-_\AA[w]$ is a Serre completely prime ideal of $U^-_\AA[w]$.

Part (2) follows from part (1), \thref{fin-ab}(4) (applied to the category $\mc{C}_w$) and 
the isomorphism $K_0( \mc{C}_w) \cong U^-_\AA[w]\spcheck$ from \thref{catUw}.
\end{proof}
It is possible that \thref{cat-Rich}(2) holds for symmetrizable Kac--Moody algebras $\g$ by arguments that 
avoid the use of \thref{KLR-cano}. 

\end{document}